# THE SPARSITY AND BIAS OF THE LASSO SELECTION IN HIGH-DIMENSIONAL LINEAR REGRESSION

By Cun-Hui Zhang[1] and Jian Huang[2]

*Rutgers University and University of Iowa*

Meinshausen and Buhlmann [*Ann. Statist.* **34** (2006) 1436–1462] showed that, for neighborhood selection in Gaussian graphical models, under a neighborhood stability condition, the LASSO is consistent, even when the number of variables is of greater order than the sample size. Zhao and Yu [(2006) *J. Machine Learning Research* **7** 2541–2567] formalized the neighborhood stability condition in the context of linear regression as a strong irrepresentable condition. That paper showed that under this condition, the LASSO selects exactly the set of nonzero regression coefficients, provided that these coefficients are bounded away from zero at a certain rate. In this paper, the regression coefficients outside an ideal model are assumed to be small, but not necessarily zero. Under a sparse Riesz condition on the correlation of design variables, we prove that the LASSO selects a model of the correct order of dimensionality, controls the bias of the selected model at a level determined by the contributions of small regression coefficients and threshold bias, and selects all coefficients of greater order than the bias of the selected model. Moreover, as a consequence of this rate consistency of the LASSO in model selection, it is proved that the sum of error squares for the mean response and the $\ell_\alpha$-loss for the regression coefficients converge at the best possible rates under the given conditions. An interesting aspect of our results is that the logarithm of the number of variables can be of the same order as the sample size for certain random dependent designs.

**1. Introduction.** Consider a linear regression model

$$(1.1) \qquad y_i = \sum_{j=1}^{p} x_{ij}\beta_j + \varepsilon_i, \qquad i = 1, \ldots, n,$$

Received January 2007; revised May 2007.
[1]Supported in part by NSF Grants DMS-05-04387 and DMS-06-04571.
[2]Supported in part by NCI/NIH Grant P30 CA 086862-06.

where $y_i$ is the response variable, $x_{ij}$ are covariates or design variables and $\varepsilon_i$ is the error term. In many applications, such as studies involving microarray or mass spectrum data, the total number of covariates $p$ can be large or even much larger than $n$, but the number of important covariates is typically smaller than $n$. With such data, regularized or penalized methods are needed to fit the model and variable selection is often the most important aspect of the analysis. The LASSO [Tibshirani (1996)] is a penalized method similar to the ridge regression but uses the $L_1$-penalty $\sum_{j=1}^{p} |\beta_j|$ instead of the $L_2$-penalty $\sum_{j=1}^{p} \beta_j^2$. An important feature of the LASSO is that it can be used for variable selection. Compared to the classical variable selection methods, such as subset selection, the LASSO has two advantages. First, the selection process in the LASSO is based on continuous trajectories of regression coefficients as functions of the penalty level and is hence more stable than subset selection methods. Second, the LASSO is computationally feasible for high-dimensional data [Osborne, Presnell and Turlach (2000a, 2000b), Efron et al. (2004)]. In contrast, computation in subset selection is combinatorial and not feasible when $p$ is large.

Several authors have studied the model-selection consistency of the LASSO in the sense of selecting exactly the set of variables with nonzero coefficients, that is, identifying the subset $\{j : \beta_j \neq 0\}$ of $\{1, \ldots, p\}$. In the low-dimensional setting with fixed $p$, Knight and Fu (2000) showed that, under appropriate conditions, the LASSO is consistent for estimating the regression parameters $\beta_j$ and their limiting distributions can have positive probability mass at 0 when $\beta_j = 0$. However, careful inspection of their results indicates that the positive probability mass at 0 is less than 1 in the limit for certain configurations of the covariates and regression coefficients, which suggests that the LASSO is not variable-selection consistent without proper assumptions. Leng, Lin and Wahba (2006) showed that the LASSO is, in general, not variable-selection consistent when the prediction accuracy is used as the criterion for choosing the penalty parameter. On the other hand, Meinshausen and Buhlmann (2006) showed that, for neighborhood selection in the Gaussian graphical models, under a neighborhood stability condition on the design matrix and certain additional regularity conditions, the LASSO is consistent, even when the number of variables tends to infinity at a rate faster than $n$. Zhao and Yu (2006) formalized the neighborhood stability condition in the context of linear regression models as a strong irrepresentable condition. They showed that under this crucial condition and certain other regularity conditions, the LASSO is consistent for variable selection, even when the number of variables $p$ is as large as $\exp(n^a)$ for some $0 < a < 1$. Thus, their results are applicable to high-dimensional regression problems, provided that the conditions, in particular, the strong irrepresentable condition, are reasonable for the data.

In this paper, we provide a different set of sufficient conditions under which the LASSO is *rate consistent* in the sparsity and bias of the selected model in high-dimensional regression. The usual definition of sparseness for model selection, as used in Meinshausen and Buhlmann (2006) and Zhao and Yu (2006), is that only a small number of regression coefficients are nonzero and all nonzero coefficients are uniformly bounded away from zero at a certain rate. Thus, variable selection is equivalent to distinguishing between nonzero and zero coefficients with a separation zone. We consider a more general concept of sparseness: a model is sparse if most coefficients are small, in the sense that the sum of their absolute values is below a certain level. Under this general sparsity assumption, it is no longer sensible to select exactly the set of nonzero coefficients. Therefore, in cases where the exact selection consistency for all $\beta_j \neq 0$ is unattainable or undesirable, we propose to evaluate the selected model with the sparsity as its dimension and the bias as the unexplained part of the mean vector and the missing large coefficients. As our goal is to select a parsimonious model which approximate the truth well, the sparsity and bias are suitable measures of performance. This is not to be confused with criteria for estimation or prediction, since we are not bound to use the LASSO for these purposes after model selection.

Under a sparse Riesz condition which limits the range of the eigenvalues of the covariance matrices of all subsets of a fixed number of covariates, we prove that the LASSO selects a model with the correct order of sparsity and controls the bias of the selected model at a level of the same order as the bias of the LASSO in the well-understood case of orthonormal design. Consequently, the LASSO selects all variables with coefficients above a threshold determined by the controlled bias of the selected model. In this sense, and in view of the optimality properties of the soft threshold method for orthonormal designs [Donoho and Johnstone (1994)], our results provide the rate consistency of the LASSO for general designs under the sparse Riesz condition. As mentioned in the previous paragraph, the LASSO does not have to be used for estimation and prediction after model selection. Nevertheless, we show that the rate consistency of the LASSO selection implies the convergence of the LASSO estimator to the true mean $Ey_i$ and coefficients $\beta_j$ at the same rate as in the case of orthonormal design.

When the number of regression coefficients exceeds the number of observations ($p > n$), there are potentially many models fitting the same data. However, there is a certain uniqueness among such models under sparsity constraints. Under the sparse Riesz condition, all sets of $q^*$ design vectors are linearly independent for a certain given rank $q^*$ so that the linear combination of design vectors is unique among all coefficient vectors of sparsity $q^*/2$ or less. Moreover, our rate consistency result proves that under mild conditions, the representation of all coefficients above a certain threshold level is determined in the selected model with high probability. Of course, such uniqueness is invalid when the sparsity assumption fails to hold.

We describe our rate consistency results in Section 2 and prove them in Section 5. Implications of the rate consistency for the convergence rate of the LASSO estimator are discussed in Section 3. The sparse Riesz and strong irrepresentable conditions do not imply each other in general, but the sparse Riesz condition is easier to interpret and less restrictive from a practical point of view. In Section 4, we provide sufficient conditions for the sparse Riesz condition for deterministic and random covariates. In Section 6, we discuss some closely related work in detail and make a few final remarks.

**2. Rate consistency of the LASSO in sparsity and bias.** The linear model (1.1) can be written as

$$(2.1) \qquad \mathbf{y} = \sum_{j=1}^{p} \beta_j \mathbf{x}_j + \boldsymbol{\varepsilon} = \mathbf{X}\boldsymbol{\beta} + \boldsymbol{\varepsilon},$$

where $\mathbf{y} \equiv (y_1, \ldots, y_n)'$, $\mathbf{x}_j$ are the columns of the design matrix $\mathbf{X} \equiv (x_{ij})_{n \times p}$, $\boldsymbol{\beta} \equiv (\beta_1, \ldots, \beta_p)'$ is the vector of regression coefficients and $\boldsymbol{\varepsilon} \equiv (\varepsilon_1, \ldots, \varepsilon_n)'$. Unless otherwise explicitly stated, we treat $\mathbf{X}$ as a given deterministic matrix.

For a given penalty level $\lambda \geq 0$, the LASSO estimator of $\boldsymbol{\beta} \in \mathbb{R}^p$ is

$$(2.2) \qquad \widehat{\boldsymbol{\beta}} \equiv \widehat{\boldsymbol{\beta}}(\lambda) \equiv \arg\min_{\boldsymbol{\beta}} \{ \|\mathbf{y} - \mathbf{X}\boldsymbol{\beta}\|^2/2 + \lambda \|\boldsymbol{\beta}\|_1 \},$$

where $\|\cdot\|$ is the Euclidean distance and $\|\boldsymbol{\beta}\|_1 \equiv \sum_j |\beta_j|$ is the $\ell_1$-norm. In this paper,

$$(2.3) \qquad \widehat{A} \equiv \widehat{A}(\lambda) \equiv \{ j \leq p : \widehat{\beta}_j \neq 0 \}$$

is considered as the model selected by the LASSO.

As mentioned in the Introduction, we consider model selection properties of the LASSO under a sparsity condition on the regression coefficients and a sparse Riesz condition on the covariates. The sparsity condition asserts the existence of an index set $A_0 \subset \{1, \ldots, p\}$ such that

$$(2.4) \qquad \#\{ j \leq p : j \notin A_0 \} = q, \qquad \sum_{j \in A_0} |\beta_j| \leq \eta_1.$$

Under this condition, there exist at most $q$ "large" coefficients and the $\ell_1$ norm of the "small" coefficients is no greater than $\eta_1$. Thus, if $q$ is of smaller order than $p$ and $\eta_1$ is small, then the high-dimensional full model $\mathbf{X}\boldsymbol{\beta}$ with $p$ coefficients can be approximated by a much lower-dimensional submodel with $q$ coefficients so that model selection makes sense. Compared with the typical assumption

$$(2.5) \qquad |A_{\boldsymbol{\beta}}| = q, \qquad A_{\boldsymbol{\beta}} \equiv \{ j : \beta_j \neq 0 \}$$

for model selection, (2.4) is mathematically weaker and much more realistic since it specifies a connected set in the parameter space $\mathbb{R}^p$ of $\boldsymbol{\beta}$. Let $(j)$ be the orderings giving $|\beta_{(1)}| \geq \cdots \geq |\beta_{(p)}|$. Another way of stating (2.4) is

$$(2.6) \qquad \sum_{j=q+1}^{p} |\beta_{(j)}| \leq \eta_1, \qquad A_0 \equiv \{(q+1), \ldots, (p)\}.$$

What should be the goal of model selection under the sparsity condition (2.4)? Unlike the usual case of (2.5), condition (2.4) allows potentially many small coefficients so that it is no longer reasonable to select exactly all variables with nonzero coefficients. Instead, a sensible goal is to select a sparse model which fits the mean vector $\mathbf{X}\boldsymbol{\beta}$ well and thus includes most (all) variables with (very) large $|\beta_j|$. Under the sparsity assumption (2.4), a natural definition of the sparsity of the selected model is $\widehat{q} = O(q)$, where

$$(2.7) \qquad \widehat{q} \equiv \widehat{q}(\lambda) \equiv |\widehat{A}| = \#\{j : \widehat{\beta}_j \neq 0\}.$$

The selected model fits the mean $\mathbf{X}\boldsymbol{\beta}$ well if its bias

$$(2.8) \qquad \widetilde{B} \equiv \widetilde{B}(\lambda) \equiv \|(\mathbf{I} - \widehat{\mathbf{P}})\mathbf{X}\boldsymbol{\beta}\|$$

is small, where $\widehat{\mathbf{P}}$ is the projection from $\mathbb{R}^n$ to the linear span of the set of selected variables $\mathbf{x}_j$ and $\mathbf{I} \equiv \mathbf{I}_n$ is the $n \times n$ identity matrix. Since the bias $\widetilde{B}$ is defined as the length of the difference between $\mathbf{X}\boldsymbol{\beta}$ and its projection to the image of $\widehat{\mathbf{P}}$, $\widetilde{B}^2$ is the sum of squares of the part of the mean vector not explained by the selected model. To measure the large coefficients for variables missing in the selected model, we define

$$(2.9) \qquad \zeta_\alpha \equiv \zeta_\alpha(\lambda) \equiv \left( \sum_{j \notin A_0} |\beta_j|^\alpha I\{\widehat{\beta}_j = 0\} \right)^{1/\alpha}, \qquad 0 \leq \alpha \leq \infty.$$

Under (2.6), $\zeta_0$ is the number of the $p$ largest $|\beta_j|$ not selected, $\zeta_2$ is the Euclidean length of these missing large coefficients and $\zeta_\infty$ is their maximum. What should be the correct order of $\widetilde{B}$ and $\zeta_\alpha$? Example 1 below indicates that under the conditions we impose, the following three quantities, or the maximum of the three, are reasonable benchmarks for $\widetilde{B}^2$ and $n\zeta_2^2$:

$$(2.10) \qquad \lambda\eta_1, \eta_2^2, \frac{q\lambda^2}{n},$$

where $\eta_2 \equiv \max_{A \subset A_0} \|\sum_{j \in A} \beta_j \mathbf{x}_j\| \leq \max_{j \leq p} \|\mathbf{x}_j\| \eta_1$.

EXAMPLE 1. Suppose we have an orthonormal design with $\mathbf{X}'\mathbf{X}/n = \mathbf{I}_p$ and i.i.d. normal error $\boldsymbol{\varepsilon} \sim N(0, \mathbf{I}_n)$. Then, (2.2) is the soft-threshold estimator [Donoho and Johnstone (1994)] with threshold level $\lambda/n$ for the individual coefficients: $\widehat{\boldsymbol{\beta}}_j = \text{sgn}(z_j)(|z_j| - \lambda/n)^+$, with $z_j \equiv \mathbf{x}'_j\mathbf{y}/n \sim N(\beta_j, 1/n)$ being the least-squares estimator of $\beta_j$. If $|\beta_j| = \lambda/n$ for $j = 1, \ldots, q + \eta_1 n/\lambda$ and $\lambda/\sqrt{n} \to \infty$, then $P\{\widehat{\beta}_j = 0\} \approx 1/2$ so that $\widetilde{B}^2 \approx 2^{-1}(q + \eta_1 n/\lambda)n(\lambda/n)^2 = 2^{-1}(q\lambda^2/n + \eta_1\lambda)$.

In this example, we observe that the order of $\widetilde{B}^2$ cannot be smaller than the first and third quantities in (2.10), while the second quantity $\eta_2^2$ is a natural choice of $\widetilde{B}^2$ as the maximum mean effect of variables with small coefficients. In the proof of Theorem 1 in Section 5 (Remark 8), we show that $\sqrt{n}\zeta_2$ is of order no greater than $\widetilde{B} + \eta_2$. Thus, we say that the LASSO is rate-consistent in model selection if, for a suitable $\alpha$ (e.g., $\alpha = 2$ or $\alpha = \infty$),

$$(2.11) \quad \widehat{q} = O(q), \qquad \widetilde{B} = O_P(B), \qquad \sqrt{n}\zeta_\alpha = O(B),$$

with the possibility of $\widetilde{B} = O(\eta_2)$ and $\zeta_\alpha = 0$ under stronger conditions, where $B \equiv \max(\sqrt{\eta_1\lambda}, \eta_2, \sqrt{q\lambda^2/n})$.

As we mentioned earlier, the main result of this paper proves the rate-consistency of the LASSO under (2.4) and a sparse Riesz condition on $\mathbf{X}$. The sparse Riesz condition controls the range of eigenvalues of covariate matrices of subsets of a fixed number of design vectors $\mathbf{x}_j$. For $A \subset \{1, \ldots, p\}$, define

$$(2.12) \quad \mathbf{X}_A \equiv (\mathbf{x}_j, j \in A), \qquad \mathbf{\Sigma}_A \equiv \mathbf{X}'_A \mathbf{X}_A / n.$$

The design matrix $\mathbf{X}$ satisfies the sparse Riesz condition (SRC) with rank $q^*$ and spectrum bounds $0 < c_* < c^* < \infty$ if

$$(2.13) \quad c_* \leq \frac{\|\mathbf{X}_A \mathbf{v}\|^2}{n\|\mathbf{v}\|^2} \leq c^* \qquad \forall A \text{ with } |A| = q^* \text{ and } \mathbf{v} \in \mathbb{R}^{q^*}.$$

Since $\|\mathbf{X}_A \mathbf{v}\|^2/n = \mathbf{v}'\mathbf{\Sigma}_A \mathbf{v}$, all the eigenvalues of $\mathbf{\Sigma}_A$ are inside the interval $[c_*, c^*]$ under (2.13) when the size of $A$ is no greater than $q^*$. While the Riesz condition asserts the equivalence of a norm $\|\sum_j v_j \xi_j\|$ and the $\ell_2$ norm $(\sum_j v_j^2)^{1/2}$ in an entire (infinite-dimensional) linear space with basis $\{\xi_1, \xi_2, \ldots\}$, the SRC provides the equivalence of the norm $\|\mathbf{\Sigma}^{1/2}\mathbf{v}\|$ and the $\ell_2$ norm $\|\mathbf{v}\|$ only in subspaces of a fixed dimension in a fixed coordinate system. The quantities $c_*$ and $c^*$ have been considered as sparse minimum and maximum eigenvalues [Meinshausen and Yu (2006), Donoho (2006)]. We call (2.13) the *sparse Riesz condition* due to its close connection to the Riesz condition as discussed above and in Section 4.2.

We prove the rate consistency (2.11) for the LASSO under the sparsity (2.4) and SRC (2.13) conditions if they are configured in certain ways between themselves and in relation to the penalty level $\lambda$. These relationships are expressed through the following ratios:

$$(2.14) \ r_1 \equiv r_1(\lambda) \equiv \left(\frac{c^*\eta_1 n}{q\lambda}\right)^{1/2}, \qquad r_2 \equiv r_2(\lambda) \equiv \left(\frac{c^*\eta_2^2 n}{q\lambda^2}\right)^{1/2}, \qquad C \equiv \frac{c^*}{c_*},$$

where $\{q, \eta_1, \eta_2, c_*, c^*\}$ are as in (2.4), (2.10) and (2.13). The quantities in (2.14) are invariant under scale changes $\{\mathbf{X}, \boldsymbol{\varepsilon}, \eta_2, \sqrt{c_*}, \sqrt{c^*}, \sqrt{\lambda}\} \to \{\mathbf{X}, \boldsymbol{\varepsilon}, \eta_2, \sqrt{c_*}, \sqrt{c^*}, \sqrt{\lambda}\}/\sigma$ and $\{\boldsymbol{\varepsilon}, \boldsymbol{\beta}, \eta_1, \eta_2, \lambda\} \to \{\boldsymbol{\beta}, \eta_1, \eta_2, \lambda\}/\sigma$. Up to the factor $c^*$ for scale adjustment, $r_1^2$ and $r_2^2$ are the ratios of the first two benchmark quantities to the third in (2.10). In terms of these scale invariant quantities, we explicitly express in our theorem the $O(1)$ in (2.11) as

(2.15) $\quad M_1^* \equiv M_1^*(\lambda) \equiv 2 + 4r_1^2 + 4\sqrt{C}r_2 + 4C,$

(2.16) $\quad M_2^* \equiv M_2^*(\lambda) \equiv \frac{8}{3}\{\frac{1}{4} + r_1^2 + r_2\sqrt{2C}(1+\sqrt{C}) + C(\frac{1}{2} + \frac{4}{3}C)\}$

and

(2.17)
$$M_3^* \equiv M_3^*(\lambda) \equiv \frac{8}{3}\Big\{\frac{1}{4} + r_1^2 + r_2\sqrt{C}(1 + 2\sqrt{1+C}) \\ + \frac{3r_2^2}{4} + C\Big(\frac{7}{6} + \frac{2}{3}C\Big)\Big\}.$$

Note that the quantities $r_j$ and $M_k^*$ in (2.14)–(2.17) are all decreasing in $\lambda$. We define a lower bound for the penalty level as

(2.18) $\qquad \lambda_* \equiv \inf\{\lambda : M_1^*(\lambda)q + 1 \leq q^*\}, \qquad \inf \varnothing \equiv \infty.$

Let $\sigma \equiv (E\|\boldsymbol{\varepsilon}\|^2/n)^{1/2}$. With the $\lambda_*$ in (2.18) and $c^*$ in (2.13), we consider the LASSO path for

(2.19) $\quad \lambda \geq \max(\lambda_*, \lambda_{n,p}), \qquad \lambda_{n,p} \equiv 2\sigma\sqrt{2(1+c_0)c^*n\log(p \vee a_n)},$

with $c_0 \geq 0$ and $a_n \geq 0$ satisfying $p/(p \vee a_n)^{1+c_0} \approx 0$. For large $p$, the lower bound here is allowed to be of the order $\lambda_{n,p} \sim \sqrt{n \log p}$ with $a_n = 0$. For example, $\lambda_* \leq \lambda_{n,p}$ if (2.13) holds for $q^* \geq (6 + 4\sqrt{C} + 4C)q + 1$, $\eta_1 \leq q\lambda_{n,p}/(nc^*)$ and $\eta_2^2 \leq q\lambda_{n,p}^2/(nc^*)$, up to $r_1 = r_2 = 1$ in (2.14). For fixed $p$, $a_n \to \infty$ is required. For i.i.d. normal errors and large $p$, the false discovery increases dramatically after the LASSO path enters the region $\lambda < \sigma\sqrt{2n \log p}$, at least in the orthonormal case.

THEOREM 1. *Let $\widehat{q}(\lambda)$, $\widetilde{B}(\lambda)$ and $\zeta_2(\lambda)$ be as in (2.7), (2.8) and (2.9), respectively, for the model $\widehat{A}(\lambda)$ selected by the LASSO with (2.2) and (2.3). Let $M_j^*$ be as in (2.15), (2.16) and (2.17). Suppose $\boldsymbol{\varepsilon} \sim N(0, \sigma^2\mathbf{I})$, $q \geq 1$, and the sparsity (2.4) and sparse Riesz (2.13) conditions hold. There then exists a set $\Omega_0$ in the sample space of $(\mathbf{X}, \boldsymbol{\varepsilon}/\sigma)$, depending on $\{\mathbf{X}\boldsymbol{\beta}, c_0, a_n\}$ only, such that*

(2.20) $\quad P\{(\mathbf{X}, \boldsymbol{\varepsilon}/\sigma) \in \Omega_0\} \geq 2 - \exp\Big(\frac{2p}{(p \vee a_n)^{1+c_0}}\Big) - \frac{2}{(p \vee a_n)^{1+c_0}} \approx 1$

*and the following assertions hold in the event* $(\mathbf{X}, \boldsymbol{\varepsilon}/\sigma) \in \Omega_0$ *for all* $\lambda$ *satisfying* (2.19):

$$(2.21) \qquad \widehat{q}(\lambda) \leq \widetilde{q}(\lambda) \equiv \#\{j : \widehat{\beta}_j(\lambda) \neq 0 \text{ or } j \notin A_0\} \leq M_1^*(\lambda)q,$$

$$(2.22) \qquad \widetilde{B}^2(\lambda) = \|(\mathbf{I} - \widehat{\mathbf{P}}(\lambda))\mathbf{X}\boldsymbol{\beta}\|^2 \leq M_2^*(\lambda)\frac{q\lambda^2}{c^*n},$$

*with* $\widehat{\mathbf{P}}(\lambda)$ *being the projection to the span of the selected design vectors* $\{\mathbf{x}_j, j \in \widehat{A}(\lambda)\}$ *and*

$$(2.23) \qquad \zeta_2^2(\lambda) = \sum_{j \notin A_0} |\beta_j|^2 I\{\widehat{\beta}_j(\lambda) = 0\} \leq M_3^*(\lambda)\frac{q\lambda^2}{c^*c_*n^2}.$$

REMARK 1. The condition $q \geq 1$ is not essential since it is only used to express quantities in Theorem 1 and its proof in terms of ratios in (2.14). Thus, (2.21), (2.22) and (2.23) are still valid for $q = 0$ if we use $r_1^2 q = c^*\eta_1 n/\lambda$ and $r_2^2 q = c^*\eta_2^2 n/\lambda^2$ to recover $M_k^* q$ from (2.15), (2.16) and (2.17), resulting in

$$\widehat{q}(\lambda) \leq 4c^*\frac{\eta_1 n}{\lambda}, \qquad \widetilde{B}^2(\lambda) \leq \frac{8}{3}\eta_1\lambda, \qquad \zeta_2^2 = 0.$$

REMARK 2. For $\eta_1 = 0$ in (2.6), we have $r_1 = r_2 = 0$ and

$$(2.24) \qquad \begin{aligned} M_1^* &= 2 + 4C, \\ M_2^* &= \frac{M_1^*}{3} + \frac{32}{9}C^2, \\ M_3^* &= \tfrac{2}{3} + \tfrac{28}{9}C + \tfrac{16}{9}C^2, \end{aligned}$$

all depend only on $C \equiv c^*/c_*$ in (2.14). In this case, (2.18) gives $\lambda_* = 0$ for $(2 + 4C)q + 1 \leq q^*$ and $\lambda_* = \infty$ otherwise. Thus, Theorem 1 requires $(2 + 4C)q + 1 \leq q^*$ in (2.4) and (2.13).

REMARK 3. The conclusions of Theorem 1 are valid for the LASSO path for all $\lambda \geq \max(\lambda_*, \lambda_{n,p})$ in the same event $(\mathbf{X}, \boldsymbol{\varepsilon}/\sigma) \in \Omega_0$. This allows data-driven selection of $\lambda$, for example, cross-validation based on prediction error. However, the theoretical justification of such a choice of $\lambda$ is unclear for model-selection purposes. Theorem 1 and simple calculation for orthonormal designs indicate that $\lambda_{n,p}$ is a good choice for model selection when $\lambda_{n,p} \geq \lambda_*$, provided we have some idea about the unknown $q$ and "known" $\{c_*, c^*, q^*\}$.

Theorem 1 is proved in Section 5. The following result is an immediate consequence of it.

THEOREM 2. *Suppose the conditions of Theorem 1 hold. Then, all variables with $\beta_j^2 > M_3^*(\lambda)q\lambda^2/\{c^*c_*n^2\}$ are selected with $j \in \widehat{A}(\lambda)$, provided $(\mathbf{X}, \varepsilon/\sigma) \in \Omega_0$ and $\lambda$ is in the interval (2.19). Consequently, if $\beta_j^2 > M_3^*(\lambda)q\lambda^2/\{c^*c_*n^2\}$ for all $j \notin A_0$, then, for all $\alpha > 0$,*

$$
(2.25) \quad \begin{aligned} P\{A_0^c \subset \widehat{A}, \widetilde{B}(\lambda) \leq \eta_2 \text{ and } \zeta_\alpha(\lambda) = 0\} \\ \geq 2 - \exp\left(\frac{2p}{(p \vee a_n)^{1+c_0}}\right) - \frac{2}{(p \vee a_n)^{1+c_0}} \approx 1. \end{aligned}
$$

Theorems 1 and 2 provide sufficient conditions under which the LASSO is rate-consistent in sparsity and bias in the sense of (2.11). It asserts that, with large probability, the LASSO selects a model with the correct order of dimension. Moreover, with large probability, the bias of the selected model is the smallest possible $\eta_2$ in the best scenario when all the large coefficients are above an explicit threshold level, and in the worst scenario, the bias is of the same order as what would be expected in the much simpler case of orthonormal design. Furthermore, with large probability, all variables with coefficients above the threshold level are selected, regardless of the values of the other coefficients. The implications of Theorem 1 on the properties of the LASSO estimator are discussed in Section 3.

In Theorems 1 and 2, conditions are imposed jointly on the design $\mathbf{X}$ and the unknown coefficients $\boldsymbol{\beta}$. Since $\mathbf{X}$ is observable, we may think of these conditions in the following way. We first impose the SRC (2.13) on $\mathbf{X}$. Given the configuration $\{q^*, c_*, c^*\}$ of the SRC and thus $C \equiv c^*/c_*$, (2.18) requires that $\{q, r_1, r_2\}$ satisfy $(2 + 4r_1^2 + 4\sqrt{C}r_2 + 4C)q + 1 \leq q^*$. Given $\{q, r_1, r_2\}$ and the penalty level $\lambda$, the condition on $\boldsymbol{\beta}$ becomes

$$|A_0^c| \leq q, \qquad \eta_1 \leq \frac{q\lambda r_1^2}{c^*n}, \qquad \eta_2^2 \leq \frac{q\lambda^2 r_2^2}{c^*n}.$$

Since Theorems 1 and 2 are valid for any fixed sample (with the exception of the "$\approx 1$" parts), $q^*, c_*, c^*, q, r_1$ and $r_2$ are all allowed to depend on $n$, but they could also be considered as fixed.

The constant factors $M_j^*$ in Theorem 1 are not sharp since crude bounds (e.g., Cauchy–Schwarz) are used several times in the proof. However, Theorem 1 is valid for any fixed $(n, p)$ with the specified configurations of the sparsity and sparse Riesz conditions. Thus, it is necessarily invariant under the scale transformations $(\mathbf{X}, \varepsilon) \to (\mathbf{X}, \varepsilon)/\sigma$ and $(\boldsymbol{\beta}', \varepsilon') \to (\boldsymbol{\beta}', \varepsilon')/\sigma$.

The SRC (2.13) is studied in Section 4 for both deterministic and random covariates. Under the Riesz condition on an infinite sequence of Gaussian covariates, we prove that (2.13) holds with fixed $0 < c_* < c^* < \infty$ and $q^* = a_0 n/\{1 \vee \log(p/n)\}$ with large probability as $(n, p) \to (\infty, \infty)$ (cf. Remark 6). This allows the application of Theorem 1 with $p$ as large as $\exp(an)$ for a small fixed $a > 0$. Section 6 contains additional discussion of our and related results after we study the LASSO estimation and SRC and prove Theorem 1.

**3. The LASSO estimation.** Here, we describe implications of Theorems 1 for the estimation properties of the LASSO. For simplicity, we confine this discussion to the special case where $c_*, c^*, r_1, r_2, c_0$ and $\sigma$ are fixed and $\lambda/\sqrt{n} \geq 2\sigma\sqrt{2(1+c_0)c^* \log p} \to \infty$. In this case, $M_k^*$ are fixed constants in (2.15), (2.16) and (2.17), and the required configurations for (2.4), (2.13) and (2.19) in Theorem 1 become

$$(3.1) \qquad M_1^* q + 1 \leq q^*, \qquad \eta_1 \leq \left(\frac{r_1^2}{c^*}\right)\frac{q\lambda}{n}, \qquad \eta_2^2 \leq \left(\frac{r_2^2}{c^*}\right)\frac{q\lambda^2}{n}.$$

Of course, $p, q$ and $q^*$ are all allowed to depend on $n$: for example, $p \gg n > q^* > q \to \infty$.

Let $A_1 \equiv \{j : \widehat{\beta}_j(\lambda) \neq 0 \text{ or } j \notin A_0\}$. Set $\mathbf{X}_1 \equiv \mathbf{X}_{A_1}$ and $\mathbf{\Sigma}_{11} \equiv \mathbf{\Sigma}_{A_1}$ as in (2.12). Define $\mathbf{b}_1 \equiv (b_j, j \in A_1)'$ for all $\mathbf{b} \in \mathbb{R}^p$. Consider the event $(\mathbf{X}, \boldsymbol{\varepsilon}/\sigma) \in \Omega_0$ in Theorem 1, in which $|A_1| \leq M_1^* q$. Since $\mathbf{\Sigma}_{11} \geq c_*$ by the SRC (2.13), the vector $\mathbf{v}_1 \equiv \mathbf{X}_1(\widehat{\boldsymbol{\beta}}_1 - \boldsymbol{\beta}_1)$ satisfies

$$(3.2) \qquad \|\mathbf{v}_1\|^2 = n\|\mathbf{\Sigma}_{11}^{1/2}(\widehat{\boldsymbol{\beta}}_1 - \boldsymbol{\beta}_1)\|^2 \geq c_* n \|\widehat{\boldsymbol{\beta}}_1 - \boldsymbol{\beta}_1\|^2.$$

The inner product of $\widehat{\boldsymbol{\beta}}_1 - \boldsymbol{\beta}_1$ and the gradient $\mathbf{g}_1 \equiv \mathbf{X}_1'(\mathbf{y} - \mathbf{X}\widehat{\boldsymbol{\beta}})$ is

$$(\widehat{\boldsymbol{\beta}}_1 - \boldsymbol{\beta}_1)'\mathbf{g}_1 = \mathbf{v}_1'(\mathbf{y} - \mathbf{X}_1\widehat{\boldsymbol{\beta}}_1) = \mathbf{v}_1'(\mathbf{X}\boldsymbol{\beta} - \mathbf{X}_1\boldsymbol{\beta}_1 + \boldsymbol{\varepsilon}) - \|\mathbf{v}_1\|^2.$$

Since $\|\mathbf{g}_1\|_\infty \leq \lambda$, and $\|\mathbf{X}\boldsymbol{\beta} - \mathbf{X}_1\boldsymbol{\beta}_1\| \leq \eta_2$,

$$(3.3) \qquad \begin{aligned} \|\mathbf{v}_1\| &\leq \|\mathbf{X}\boldsymbol{\beta} - \mathbf{X}_1\boldsymbol{\beta}_1 + \mathbf{P}_1\boldsymbol{\varepsilon}\| + n^{-1/2}\|\mathbf{\Sigma}_{11}^{-1/2}\mathbf{g}_1\| \\ &\leq \eta_2 + \|\mathbf{P}_1\boldsymbol{\varepsilon}\| + \lambda\left(\frac{|A_1|}{c_* n}\right)^{1/2}, \end{aligned}$$

where $\mathbf{P}_1 \equiv \mathbf{X}_1'\mathbf{\Sigma}_{11}^{-1}\mathbf{X}_1/n$ is the projection to the range of $\mathbf{X}_1$. Since $\text{rank}(\mathbf{P}_1) = |A_1| \leq M_1^* q$, we are able to show that $\|\mathbf{P}_1\boldsymbol{\varepsilon}\|$ is of the order $\sigma\sqrt{q \log p}$ under the normality assumption. Thus, (3.2) and (3.3) lead to Theorem 3 below. The inequality (2.21) plays a crucial role here since it controls $|A_1|$ and then allows the application of the SRC.

THEOREM 3. *Let $c_*, c^*, r_1, r_2, c_0$ and $\sigma$ be fixed and $1 \leq q \leq p \to \infty$. Let $\lambda = 2\sigma\sqrt{2(1+c_0')c^* n \log p}$ with a fixed $c_0' \geq c_0$ and $\Omega_0$ be as in Theorem 1. Suppose the conditions of Theorem 1 hold with configurations satisfying* (3.1). There then exist constants $M_k^*$ depending only on $c_*, c^*, r_1, r_2$ and $c_0'$ and a set $\widetilde{\Omega}_q$ in the sample space of $(\mathbf{X}, \boldsymbol{\varepsilon}/\sigma)$ depending only on $q$ such that

$$
(3.4) \quad P\{(\mathbf{X}, \boldsymbol{\varepsilon}/\sigma) \notin \Omega_0 \cap \widetilde{\Omega}_q | \mathbf{X}\}
$$
$$
\leq e^{2/p^{c_0}} - 1 + \frac{2}{p^{1+c_0}} + \left(\frac{1}{p^2} + \frac{\log p}{p^2/4}\right)^{(q+1)/2} \to 0
$$

and the following assertions hold in the event $(\mathbf{X}, \boldsymbol{\varepsilon}/\sigma) \in \Omega_0 \cap \widetilde{\Omega}_q$:

$$
(3.5) \quad \|\mathbf{X}(\widehat{\boldsymbol{\beta}} - \boldsymbol{\beta})\| \leq M_4^* \sigma \sqrt{q \log p}
$$

and, for all $\alpha \geq 1$,

$$
(3.6) \quad \|\widehat{\boldsymbol{\beta}} - \boldsymbol{\beta}\|_\alpha \equiv \left(\sum_{j=1}^p |\widehat{\beta}_j - \beta_j|^\alpha\right)^{1/\alpha} \leq M_5^* \sigma q^{1/(\alpha \wedge 2)} \sqrt{(\log p)/n}.
$$

REMARK 4. The convergence rates in (3.5) and (3.6) are sharp for the LASSO under the given conditions since the convergence rate for (3.6) is $q^{1/\alpha}(\lambda/n + \sigma/\sqrt{n})$, $1 \leq \alpha \leq 2$, for orthogonal designs and the bias for a single $\widehat{\beta}_j$ could be of the order $\sqrt{q(\log p)/n}$, even under the strong irrepresentable condition. Moreover, by Foster and George (1994), the risk inflation factor $\sqrt{\log p}$ is optimal for (3.5) and (3.6) with $\alpha = 2$. We discuss related work in Section 6 after we study the SRC and prove Theorem 1.

PROOF OF THEOREM 3. Define $\mathbf{P}_A \equiv \mathbf{X}_A' \boldsymbol{\Sigma}_A^{-1} \mathbf{X}_A / n$ with the notation in (2.12) and

$$
\widetilde{\Omega}_q \equiv \left\{\max_{q < |A| \leq p} \frac{\|\mathbf{P}_A \boldsymbol{\varepsilon}\|^2}{\sigma^2 |A|} \leq 4 \log p\right\}.
$$

For deterministic $A$ with $\text{rank}(\mathbf{X}_A) = m$, $\|\mathbf{P}_A \boldsymbol{\varepsilon}\|^2/\sigma^2 \sim \chi_m^2$ so that

$$
P\{\|\mathbf{P}_A \boldsymbol{\varepsilon}\|^2/\sigma^2 \geq m(1 + 4\log p)\} \leq \{p^{-4}(1 + 4\log p)\}^{m/2},
$$

by the standard large deviation inequality. It follows that

$$
1 - P\{\widetilde{\Omega}_q\} \leq \sum_{m=q+1}^p \binom{p}{m} \{p^{-4}(1 + 4\log p)\}^{m/2} \leq \left(\frac{1}{p^2} + \frac{\log p}{p^2/4}\right)^{(q+1)/2},
$$

due to the facts that $\binom{p}{m} \leq p^m/m!$ and $1 + 4\log p \leq p^2$. Since $q + 1 \leq q^*$, the arguments for (3.2) and (3.3) are still valid if we require $|A_1| \geq q + 1$ (making $A_1$ larger). Thus, (3.5) follows from (2.21) and (3.3), due to $\|\mathbf{P}_1 \boldsymbol{\varepsilon}\| \leq 2\sigma\sqrt{|A_1|\log p}$ in $\widetilde{\Omega}_q$. Similarly, by both (3.2) and (3.3), we have, in $\Omega_0 \cap \widetilde{\Omega}_q$,

$$
\left(\sum_{j \in A_1} |\widehat{\beta}_j - \beta_j|^2\right)^{1/2} \leq O(1) \sigma \sqrt{|A_1|(\log p)/n}
$$

uniformly. Thus, since $A_1^c \subseteq A_0$, (3.6) follows from

$$(3.7) \qquad \left(\sum_{j \in A_0} |\beta_j|^\alpha\right)^{1/\alpha} \leq O(1)\sigma q^{1/(\alpha \wedge 2)}\sqrt{(\log p)/n}$$

for $\alpha = 1, 2$ and $\alpha = \infty$. For $\alpha = 1$, (3.7) follows from the second inequality of (3.1). For $\alpha = 2$, $\#\{j \in A_0 : |\beta_j| > \lambda/n\} = O(q)$, by (3.7) for $\alpha = 1$, so that, by the SRC (2.13) and the third inequality of (3.1),

$$\sum_{j \in A_0} \beta_j^2 I\{|\beta_j| > \lambda/n\} \leq O(1/n)\left\|\sum_{j \in A_0} \beta_j \mathbf{x}_j I\{|\beta_j| > \lambda/n\}\right\|^2$$

$$\leq O(\eta_2^2/n) = O(q\lambda^2/n^2).$$

Thus, (3.7) for $\alpha = 2$ follows from $\alpha = 1$. Finally, (3.7) for $\alpha = \infty$ follows from $\beta_j^2 \leq \|\beta_j \mathbf{x}_j\|^2/(nc_*) \leq \eta_2^2/(nc_*)$. □

**4. The sparse Riesz condition.** In this section, we provide sufficient conditions for the sparse Riesz condition. We divide the section into two subsections respectively for deterministic and random design matrices $\mathbf{X}$. In the case of random design, the rows of $\mathbf{X}$ are assumed to be i.i.d. vectors, but the entries within a row are allowed to be dependent.

We consider the sparse Riesz condition (2.13) and its general version

$$(4.1) \quad c_*(m) \equiv \min_{|A|=m} \min_{\|\mathbf{v}\|=1} \|\mathbf{X}_A \mathbf{v}\|^2/n, \ c^*(m) \equiv \max_{|A|=m} \max_{\|\mathbf{v}\|=1} \|\mathbf{X}_A \mathbf{v}\|^2/n,$$

for ranks $0 \leq m \leq p$, with the convention that $c_*(0) \equiv c^*(0) \equiv \sqrt{c_*(1)c^*(1)}$. This includes (2.13) with $c_* = c_*(q^*)$ and $c^* = c^*(q^*)$. As we mentioned earlier, (4.1) reduces to the requirement that all of the eigenvalues of $\mathbf{\Sigma}_A$ in (2.12) lie in the interval $[c_*(m), c^*(m)]$ when $|A| \leq m$. If $\mathbf{x}_j$ are standardized with $\mathbf{x}_j' \mathbf{x}_j/n = 1$, then $c_*(1) = c^*(1) = 1$. In general, $c_*(1) \leq \|\mathbf{x}_j\|^2/n \leq c^*(1)$. It is clear that $c_*(m)$ is decreasing in $m$ with $c_*(n+1) = 0$, $c^*(m)$ is increasing in $m$ and the Cauchy-Schwarz inequality gives the subadditivity $c^*(m_1 + m_2) \leq c^*(m_1) + c^*(m_2)$.

4.1. *Deterministic design matrices.* Proposition 1 below provides a simple sufficient condition for (2.13). It is actually an $\ell_\alpha$-version of Geršgorin's theorem.

PROPOSITION 1. *Suppose that* $\mathbf{X}$ *is standardized with* $\|\mathbf{x}_j\|^2/n = 1$. *Let* $\rho_{jk} = \mathbf{x}_j' \mathbf{x}_k/n$ *be the correlation. If*

$$(4.2) \qquad \max_{|A|=q^*} \inf_{\alpha \geq 1} \left\{\sum_{j \in A}\left(\sum_{k \in A, k \neq j} |\rho_{jk}|^{\alpha/(\alpha-1)}\right)^{\alpha-1}\right\}^{1/\alpha} \leq \delta < 1,$$

*then the sparse Riesz condition (2.13) holds with rank $q_*$ and spectrum bounds $c_* = 1 - \delta$ and $c^* = 1 + \delta$. In particular, (2.13) holds with $c_* = 1 - \delta$ and $c^* = 1 + \delta$ if*

$$\text{(4.3)} \qquad \max_{1 \le j < k \le p} |\rho_{jk}| \le \frac{\delta}{q^* - 1}, \qquad \delta < 1.$$

REMARK 5. If $\delta = 1/3$, then $C \equiv c^*/c_* = 2$ and Theorem 1 is applicable for (1.1) if $10q + 1 \le q^*$ and $\eta_1 = 0$ in (2.4).

PROOF OF PROPOSITION 1. Let $\boldsymbol{\Sigma}_A = (\rho_{jk})_{j \in A, k \in A}$ be the covariance matrix for variables in $A$, as in (2.12). Let $|A| = q^*$ and $\mathbf{b} = (b_1, \ldots, b_{q^*})$ be an eigenvector of $\boldsymbol{\Sigma}_A$ with eigenvalue $\tau$. Then,

$$b_j + \sum_{k \ne j, k \in A} \rho_{jk} b_k = \tau b_j$$

so that, by the Hölder inequality,

$$\sum_{j \in A} |(1 - \tau) b_j|^\alpha = \sum_{j \in A} \left| \sum_{k \ne j} \rho_{jk} b_k \right|^\alpha \le \sum_{j \in A} \left( \sum_{k \ne j} |\rho_{jk}|^{\alpha/(\alpha-1)} \right)^{\alpha - 1} \sum_{k \in A} |b_k|^\alpha.$$

After the cancellation of $\sum_{k \in A} |b_k|^\alpha$, we find, by (4.2), that $|1 - \tau| \le \delta$. This gives (2.13) with $c_* = 1 - \delta$ and $c^* = 1 + \delta$ as the interval $[c_*, c^*]$ contains all eigenvalues of $\boldsymbol{\Sigma}_A$ with $|A| = q^*$. If (4.3) holds, then, as $\alpha \to \infty$,

$$\left\{ \sum_{j \in A} \left( \sum_{k \in A, k \ne j} |\rho_{jk}|^{\alpha/(\alpha-1)} \right)^{\alpha - 1} \right\}^{1/\alpha} \le \frac{\delta}{q^* - 1} \{q^*(q^* - 1)^{\alpha - 1}\}^{1/\alpha}$$

$$= \delta (q^*)^{1/\alpha} (q^* - 1)^{-1/\alpha} \to \delta.$$

The proof of Proposition 1 is complete. $\square$

4.2. *Random design matrices.* Suppose we would like to investigate the linear relationships between a response variable $Y$ and infinitely many possible covariates $\{\xi_k, k = 1, 2, \ldots\}$. Suppose that in the $n$th experiment, we collect a sample from the dependent variable $Y$ and $p$ covariates so that we observe $n$ independent copies $(y^{(n)}, x_{ij}^{(n)}, j = 1, \ldots, p^{(n)})$ of the random vector $(Y, \xi_{k_j}, j = 1, \ldots, p)$ for certain $k_1 < \cdots < k_p$, $p \equiv p^{(n)}$. In this case, the linear model (1.1) becomes

$$\text{(4.4)} \qquad y_i^{(n)} = \sum_{j=1}^p \beta_j^{(n)} x_{ij}^{(n)} + \varepsilon_i^{(n)}.$$

In what follows, the superscript $^{(n)}$ is often omitted.

The infinite population sequence $\{\xi_j, j = 1, 2, \ldots\}$ satisfies the Riesz condition if there exist fixed $0 < \rho_* < \rho^* < \infty$ such that

$$(4.5) \qquad \rho_* \sum_{j=1}^{\infty} b_j^2 \leq E \left| \sum_{j=1}^{\infty} b_j \xi_j \right|^2 \leq \rho^* \sum_{j=1}^{\infty} b_j^2$$

for all constants $b_j$. Let $\mathbf{x}^i \equiv (x_{i1}^{(n)}, \ldots, x_{ip}^{(n)})$ be the row vectors of $\mathbf{X} \equiv (x_{ij}^{(n)})_{n \times p} = (\mathbf{x}_1, \ldots, \mathbf{x}_p)$ in (4.4). Since $\mathbf{x}^i, i = 1, \ldots, n$, are i.i.d. copies of $(\xi_{k_1}, \ldots, \xi_{k_p})$, (4.5) implies that

$$\rho_* \|\mathbf{b}\|^2 \leq E \left| \sum_{j=1}^{p} b_j \xi_{k_j} \right|^2 = E \sum_{i=1}^{n} \frac{(\mathbf{b}' \mathbf{x}^i)^2}{n} = \frac{E \|\mathbf{X} \mathbf{b}\|^2}{n} \leq \rho^* \|\mathbf{b}\|^2.$$

However, this does not guarantee that $0 < \kappa \leq c_*(m) < c^*(m) \leq 1/\kappa$ with large probability for all $m$. In particular, we always have $c_*(n+1) = 0$.

PROPOSITION 2. *Suppose that the $n$ rows of a random matrix $\mathbf{X}_{n \times p}$ are i.i.d. copies of a subvector $(\xi_{k1}, \ldots, \xi_{k_p})$ of a zero-mean random sequence $\{\xi_j, j = 1, 2, \ldots\}$ satisfying (4.5). Let $c_*(m)$ and $c^*(m)$ be as in (4.1).*

(i) *Suppose $\{\xi_k, k \geq 1\}$ is a Gaussian sequence. Let $\epsilon_k$, $k = 1, 2, 3, 4$, be positive constants in $(0, 1)$ satisfying $m \leq \min(p, \epsilon_1^2 n)$, $\epsilon_1 + \epsilon_2 < 1$ and $\epsilon_3 + \epsilon_4 = \epsilon_2^2/2$. Then, for all $(m, n, p)$ satisfying $\log \binom{p}{m} \leq \epsilon_3 n$,*

$$(4.6) \qquad P\{\tau_* \rho_* \leq c_*(m) \leq c^*(m) \leq \tau^* \rho^*\} \geq 1 - 2 e^{-n \epsilon_4},$$

*where $\tau_* \equiv (1 - \epsilon_1 - \epsilon_2)^2$ and $\tau^* \equiv (1 + \epsilon_1 + \epsilon_2)^2$.*

(ii) *Suppose $\max_{j \leq p} \|\xi_{k_j}\|_\infty \leq K_n < \infty$. Then, for any $\tau_* < 1 < \tau^*$, there exists a constant $\epsilon_0 > 0$ depending only on $\rho^*, \rho_*, \tau_*$ and $\tau^*$ such that*

$$P\{\tau_* \rho_* \leq c_*(m) \leq c^*(m) \leq \tau^* \rho^*\} \to 1$$

*for $m \equiv m_n \leq \epsilon_0 K_n^{-1} \sqrt{n/\log p}$, provided $\sqrt{n}/K_n \to \infty$.*

REMARK 6. By the Stirling formula, for $p/n \to \infty$,

$$m \leq \epsilon_3 n / \log(p/n) \Rightarrow \log \binom{p}{m} \leq (\epsilon_3 + o(1)) n.$$

Thus, Proposition 2(i) is applicable up to $p = e^{an}$ for some small $a > 0$.

REMARK 7. Supposing $m = p$, $p/n \to \epsilon_1^2 \in (0, 1)$ and $\xi_j$ are i.i.d. $N(0, 1)$, Geman (1980) proved $c^*(m) \to (1 + \epsilon_1)^2$ and Silverstein (1985) proved $\rho_* \to (1 - \epsilon_1)^2$. Silverstein's results can be directly used to prove bounds similar to (4.6) [cf. Zhang and Huang (2006)]. We refer to Bai (1999) and Davidson and Szarek (2001) for further discussion on random covariance matrices.

PROOF OF PROPOSITION 2. (i) Let $S^{m-1}$ be the unit sphere of $\mathbb{R}^m$ and $\mathbf{P}_m : \mathbb{R}^p \to \mathbb{R}^m$ be $m \times p$ projection matrices taking $m$ out of $p$ coordinates of $\mathbb{R}^p$. Define

$$\tau_-(\mathbf{P}_m) \equiv \inf_{\mathbf{b} \in S^{m-1}} \frac{\|\mathbf{X}\mathbf{P}'_m \mathbf{b}\|^2}{E\|\mathbf{X}\mathbf{P}'_m \mathbf{b}\|^2}, \qquad \tau_+(\mathbf{P}_m) \equiv \sup_{\mathbf{b} \in S^{m-1}} \frac{\|\mathbf{X}\mathbf{P}'_m \mathbf{b}\|^2}{E\|\mathbf{X}\mathbf{P}'_m \mathbf{b}\|^2}.$$

Since $\rho_* \leq E\|\mathbf{X}\mathbf{P}'_m \mathbf{b}\|^2/n \leq \rho^*$, by (4.1), we have

$$\begin{aligned}
P_{m,n,p} &\equiv P\{\tau_* \rho_* \leq c_*(m) \leq c^*(m) \leq \tau^* \rho^*\} \\
&\geq P\Big\{\tau_* \leq \min_{\mathbf{P}_m} \tau_-(\mathbf{P}_m) \leq \max_{\mathbf{P}_m} \tau_+(\mathbf{P}_m) \leq \tau^*\Big\}.
\end{aligned} \tag{4.7}$$

For a fixed $\mathbf{P}_m$, let $\mathbf{\Sigma}_m$ be the $m \times m$ population covariance matrices of the rows of $\mathbf{X}\mathbf{P}'_m$ and $\mathbf{U} \equiv \mathbf{X}\mathbf{P}'_m \mathbf{\Sigma}_m^{-1/2}$. Since $\mathbf{U}$ is then an $n \times m$ matrix of $N(0,1)$,

$$\tau_+(\mathbf{P}_m) = \sup_{\mathbf{b} \in S^{m-1}} \frac{\|\mathbf{U}\mathbf{\Sigma}_m^{1/2}\mathbf{b}\|^2}{n\|\mathbf{\Sigma}_m^{1/2}\mathbf{b}\|^2} = \sup_{\mathbf{b} \in S^{m-1}} \|\mathbf{U}\mathbf{b}\|^2/n = \lambda_{\max}(\mathbf{W}/n)$$

and $\tau_-(\mathbf{P}_m) = \lambda_{\min}(\mathbf{W}/n)$, where $\mathbf{W} \equiv \mathbf{U}'\mathbf{U}$ is an $m \times m$ matrix with the Wishart distribution $W_m(\mathbf{I}, n)$ [cf. Eaton (1983)]. Since $m/n \leq \epsilon_1^2$, for the prescribed $\tau_*$ and $\tau^*$, Theorem II.13 of Davidson and Szarek (2001) gives

$$\max(P\{\lambda_{\min}(\mathbf{W}/n) \leq \tau_*\}, P\{\lambda_{\max}(\mathbf{W}/n) \geq \tau^*\}) \leq e^{-n\epsilon_2^2/2}.$$

Thus, since there exist a total of $\binom{p}{m}$ choices of $\mathbf{P}_m$, by (4.7),

$$\begin{aligned}
-P_{m,n,p} &\leq \binom{p}{m}(1 - P\{\tau_* \leq \lambda_{\min}(\mathbf{W}/n) \leq \lambda_{\max}(\mathbf{W}/n) \leq \tau^*\}) \\
&\leq 2\binom{p}{m} e^{-n\epsilon_2^2/2} \leq 2e^{-n\epsilon_4}.
\end{aligned} \tag{4.8}$$

(ii) Define $f_n(\mathbf{b}) \equiv (\|\mathbf{X}\mathbf{P}'_m \mathbf{b}\|^2/n)^{1/2}$ and $f(\mathbf{b}) \equiv (Ef_n^2(\mathbf{b}))^{1/2}$. By (4.5), $f^2(\mathbf{b})/\|\mathbf{b}\|^2 \in [\rho_*, \rho^*]$ for all $\mathbf{b} \neq 0$. Since both $f_n$ and $f$ are norms in $\mathbb{R}^m$,

$$\left|\frac{f_n(\mathbf{b} + \widetilde{\mathbf{b}})}{f(\mathbf{b} + \widetilde{\mathbf{b}})} - \frac{f_n(\mathbf{b})}{f(\mathbf{b})}\right| \leq \left(\frac{f_n(\widetilde{\mathbf{b}})}{f(\widetilde{\mathbf{b}})} + \frac{f_n(\mathbf{b})}{f(\mathbf{b})}\right)\frac{f(\widetilde{\mathbf{b}})}{f(\mathbf{b} + \widetilde{\mathbf{b}})}.$$

Let $S_{\epsilon_1}^{m-1}$ be an $\epsilon_1$-net in $S^{m-1}$ with $2\epsilon_1\sqrt{\rho^*/\rho_*} \leq 1/5$. We have

$$\tau_+^{1/2}(\mathbf{P}_m) \leq \max_{\mathbf{b} \in S_{\epsilon_1}^{m-1}} \frac{f_n(\mathbf{b})}{f(\mathbf{b})} + 2\tau_+^{1/2}(\mathbf{P}_m)\epsilon_1\sqrt{\frac{\rho^*}{\rho_*}} \leq \frac{5}{4}\max_{\mathbf{b} \in S_{\epsilon_1}^{m-1}} \frac{f_n(\mathbf{b})}{f(\mathbf{b})}$$

and

$$\tau_-^{1/2}(\mathbf{P}_m) \geq \min_{\mathbf{b} \in S_{\epsilon_1}^{m-1}} \frac{f_n(\mathbf{b})}{f(\mathbf{b})} - \frac{1}{5}\tau_+^{1/2}(\mathbf{P}_m).$$

Since $f_n^2(\mathbf{b})/f^2(\mathbf{b})$ is the average of $n$ i.i.d. variables, each with mean 1 and uniformly bounded by $mK_n^2/\rho_*$, by the Bernstein inequality, we have

$$P\{|f_n^2(\mathbf{b})/f^2(\mathbf{b}) - 1| > 7/25\} \leq 2\exp\left(\frac{-\epsilon_2 n}{mK_n^2}\right)$$

for certain $\epsilon_2$ depending on $\rho_*$ only. Thus, for $\tau^* = (5/4)^2(1 + 7/25) = 2$ and $\tau_* = (\sqrt{1 - 7/25} - \sqrt{2}/5)^2 = 8/25$, we have

$$1 - P_{m,n,p} \leq 2\binom{p}{m}|S_{\epsilon_1}^{m-1}|\exp\left(\frac{-\epsilon_2 n}{mK_n^2}\right).$$

Since $|S_{\epsilon_1}^{m-1}|/m! = O(1)$, $P_{m,n,p} \to 1$ for $\epsilon_2 n/(mK_n^2) > 2m\log p$. This proves (ii) for the specific $\{\tau_*, \tau^*\}$. We omit the proof for the general $\{\tau_*, \tau^*\}$. □

**5. Proof of Theorem 1.** Taking the scale change $\{\boldsymbol{\varepsilon}, \boldsymbol{\beta}, \lambda\} \to \{\boldsymbol{\varepsilon}/\sigma, \boldsymbol{\beta}/\sigma, \lambda/\sigma\}$ if necessary, we assume $\boldsymbol{\varepsilon} \sim N(0, \mathbf{I})$, without loss of generality. It follows from the Karush–Kuhn–Tucker condition that a vector $\mathbf{b} \equiv (b_1, \ldots, b_p)'$ is the solution $\widehat{\boldsymbol{\beta}}$ of (2.2) and only if

$$(5.1) \quad \begin{cases} \mathbf{x}_j'(\mathbf{y} - \mathbf{X}\mathbf{b}) = \operatorname{sgn}(b_j)\lambda, & |b_j| > 0, \\ |\mathbf{x}_j'(\mathbf{y} - \mathbf{X}\mathbf{b})| \leq \lambda, & b_j = 0. \end{cases}$$

This allows us to define slightly more general versions of the $\widehat{A}$ in (2.3) and its dimension as

$$(5.2) \quad \{j : \widehat{\beta}_j \neq 0\} \subseteq A_1 \subseteq \{j : |\mathbf{x}_j'(\mathbf{y} - \mathbf{X}\widehat{\boldsymbol{\beta}})| = \lambda\} \cup A_0^c, \qquad q_1 \equiv |A_1|.$$

Set $A_2 \equiv \{1, \ldots, p\} \setminus A_1$, $A_3 \equiv A_1 \setminus A_0$, $A_4 \equiv A_1 \cap A_0$, $A_5 \equiv A_2 \setminus A_0$ and $A_6 \equiv A_2 \cap A_0$. For $A_k \subset A_j$, let $\mathbf{Q}_{kj}$ be the matrix representing the selection of variables in $A_k$ from $A_j$, defined as $\mathbf{Q}_{kj}\boldsymbol{\beta}_j = \boldsymbol{\beta}_k$, where $\boldsymbol{\beta}_k \equiv (\beta_j, j \in A_k)$. For example, $\boldsymbol{\beta}_1' = \boldsymbol{\beta}_3'\mathbf{Q}_{31} + \boldsymbol{\beta}_4'\mathbf{Q}_{41}$ since $A_1 = A_3 \cup A_4$ and $A_3 \cap A_4 = \varnothing$. We define matrices $\boldsymbol{\Sigma}_{jk} \equiv n^{-1}\mathbf{X}_j'\mathbf{X}_k$, and the projection $\mathbf{P}_1$ from $\mathbb{R}^n$ to the span of $\{\mathbf{x}_j, j \in A_1\}$. We apply all arithmetic and logic operations and univariate functions to vectors componentwise. For example, $\mathbf{v} \times |\boldsymbol{\beta}| = (v_1|\beta_1|, \ldots, v_p|\beta_p|)'$. The SRC (4.1) for a general rank $m$ is used in most parts of the proof, rather that (2.13). Table 1 summarizes the meanings of the index sets $A_j$.

We note that $\widehat{q} = q_1$ and $\widehat{\mathbf{P}} = \mathbf{P}_1$ when we choose the smallest possible $A_1$ in (5.2) and that $A_5 = \varnothing$ when we choose the largest possible $A_1$. In our analysis of the LASSO, quantities related to the coefficients in the sets $A_j, j = 0, 1, 2$, are often decomposed into those involving the more specific sets $A_j, j = 3, 4, 5, 6$.

It follows from (5.1) that

$$(5.3) \quad \mathbf{s}_j \equiv \mathbf{X}_{A_j}'(\mathbf{y} - \mathbf{X}\widehat{\boldsymbol{\beta}})/\lambda \in [-1, 1], \qquad j = 1, 3, 4.$$

TABLE 1
*Sets of variables considered in the proof*

|  | "Large" $\|\beta_j\|$ <br> $j \notin A_0$ | "Small" $\|\beta_j\|$ <br> $j \in A_0$ | **Quantities** <br> **to be bounded** |
|---|---|---|---|
| $A_1$ : selected $j$ and some $j \notin A_0$ | $A_3$ | $A_4$ | $\widehat{q} \leq q_1 \equiv \|A_1\|$ |
| $A_2$ : $j$ not in $A_1$ | $A_5$ | $A_6$ | $\|(\mathbf{I} - \widehat{\mathbf{P}})\mathbf{X}\boldsymbol{\beta}\|$ |

Our goal is to find upper bounds for the dimension $q_1 \equiv |A_1|$ and the bias terms $\|(\mathbf{I} - \mathbf{P}_1)\mathbf{X}\boldsymbol{\beta}\|$ and $\|\boldsymbol{\beta}_5\|$ for all the $A_1$ in (5.2). By (5.1), (5.2) and Table 1, we have $|\mathbf{s}_4| = 1$ for each component so that $\|\mathbf{s}_4\|^2 = |A_4|$ and $q_1 \equiv |A_1| = |A_3| + |A_4| \leq q + \|\mathbf{s}_4\|^2$. Our plan is to find upper bounds for the lengths of the vectors $\mathbf{v}_{14}$, $\mathbf{w}_2$ and $\boldsymbol{\beta}_5$, where

$$(5.4) \qquad \mathbf{v}_{1j} \equiv \frac{\lambda}{n^{1/2}} \boldsymbol{\Sigma}_{11}^{-1/2} \mathbf{Q}_{j1}' \mathbf{s}_j, \qquad \mathbf{w}_k \equiv (\mathbf{I} - \mathbf{P}_1)\mathbf{X}_k \boldsymbol{\beta}_k,$$

for $j = 3, 4$ and $k = 2, \ldots, 6$. Since $\mathbf{X}\boldsymbol{\beta} = \mathbf{X}_1\boldsymbol{\beta}_1 + \mathbf{X}_2\boldsymbol{\beta}_2$ and $(\mathbf{I} - \mathbf{P}_1)\mathbf{X}_1\boldsymbol{\beta}_1 = 0$, by (5.4) and (4.1), the fact that $\|\mathbf{s}_4\|^2 = |A_4|$ implies that

$$(5.5) \qquad \|\mathbf{v}_{14}\|^2 \geq \frac{\lambda^2(q_1 - q)}{nc^*(q_1)}, \qquad \|\mathbf{w}_2\|^2 = \|(\mathbf{I} - \mathbf{P}_1)\mathbf{X}\boldsymbol{\beta}\|^2.$$

Thus, we proceed to find upper bounds for $\|\mathbf{v}_{14}\|$, $\|\mathbf{w}_2\|$ and $\|\boldsymbol{\beta}_5\|$.

We divide the rest of the proof into three steps. Step 1 proves that the quadratic $\|\mathbf{v}_{14}\|^2 + \|\mathbf{w}_2\|^2$ is no greater than a linear function of $\{\|\mathbf{v}_{14}\|, \|\mathbf{w}_2\|, \|\boldsymbol{\beta}_5\|_1, \|\mathbf{P}_1\mathbf{X}_2\boldsymbol{\beta}_2\|\}$ with a stochastic slope. This step is crucial since the identity and inequalities in the Karush–Kuhn–Tucker (5.1) must be combined in a proper way to cancel out the cross-product term of $\mathbf{s}_4$ and $\boldsymbol{\beta}_5$. Step 2 translates the results of Step 1 into upper bounds for $q_1$, $\|\mathbf{w}_2\|^2$ and $\|\boldsymbol{\beta}_5\|^2$, essentially with careful applications of the Cauchy–Schwarz inequality, for a suitable level of the random slope and the prescribed penalty levels $\lambda$. The upper bounds in Step 2 are of the same form as in the conclusions of the theorem, but still involve $c_*(|A|)$ and $c^*(|A|)$ with random $A \subset A_1 \cup A_5$ instead of the $c_*$ and $c^*$ specified in (2.13). Step 3 completes the proof by finding probabilistic bounds for the random slope and by showing $|A_1 \cup A_5| \leq q^*$ for the rank $q^*$ in (2.13). We need a lemma for the interpretation of (4.1).

LEMMA 1. *Let $c_*(m)$ and $c^*(m)$ be as in (4.1). Let $A_k \subset \{1, \ldots, p\}$, $\mathbf{X}_k \equiv (\mathbf{x}_j, j \in A_k)$ and $\boldsymbol{\Sigma}_{1k} \equiv \mathbf{X}_1' \mathbf{X}_k / n$. Then,*

$$(5.6) \quad \frac{\|\mathbf{v}\|^2}{c^*(|A_1|)} \leq \|\boldsymbol{\Sigma}_{11}^{-1/2}\mathbf{v}\|^2 \leq \frac{\|\mathbf{v}\|^2}{c_*(|A_1|)}, \qquad \|\boldsymbol{\beta}_k\|_1^2 \leq \frac{\|\mathbf{X}_k\boldsymbol{\beta}_k\|^2 |A_k|}{nc_*(|A_k|)},$$

*for all* $\mathbf{v}$ *of proper dimension. Furthermore, if* $A_k \cap A_1 = \varnothing$, *then*

$$\|\boldsymbol{\beta}_k\|^2 + \|\boldsymbol{\Sigma}_{11}^{-1}\boldsymbol{\Sigma}_{1k}\boldsymbol{\beta}_k\|^2 \leq \frac{\|(\mathbf{I} - \mathbf{P}_1)\mathbf{X}_k\boldsymbol{\beta}_k\|^2}{nc_*(|A_1 \cup A_k|)}, \quad (5.7)$$

*where* $\mathbf{P}_1$ *is the projection to the span of* $\{\mathbf{x}_j, j \in A_1\}$.

REMARK 8. For $A_5 \equiv \{j : j \notin A_0, \widehat{\beta}_j = 0\}$, Lemma 1 gives $\zeta_2^2 = \|\boldsymbol{\beta}_5\|^2 \leq (\widetilde{B} + \eta_2)^2/(nc_*)$, provided $|A_1 \cup A_5| \leq q^*$ under the SRC (2.13).

PROOF OF LEMMA 1. We only prove the inequality of (5.7), since the rest of the lemma follows directly from the Cauchy–Schwarz inequality and (4.1). Let $\mathbf{v} \equiv -\boldsymbol{\Sigma}_{11}^{-1}\boldsymbol{\Sigma}_{1k}\boldsymbol{\beta}_k$. Since $(\mathbf{I} - \mathbf{P}_1)\mathbf{X}_k\boldsymbol{\beta}_k = \mathbf{X}_1\mathbf{v} + \mathbf{X}_k\boldsymbol{\beta}_k$,

$$\|(\mathbf{I} - \mathbf{P}_1)\mathbf{X}_k\boldsymbol{\beta}_k\|^2 = (\mathbf{v}', \boldsymbol{\beta}_k')(\mathbf{X}_1, \mathbf{X}_k)'(\mathbf{X}_1, \mathbf{X}_k)\begin{pmatrix}\mathbf{v}\\\boldsymbol{\beta}_k\end{pmatrix}$$

$$\geq nc_*(|A_1 \cup A_k|)(\|\mathbf{v}\|^2 + \|\boldsymbol{\beta}_k\|^2).$$

The proof of Lemma 1 is complete. □

*Step* 1. In this step, we prove

$$\|\mathbf{v}_{14}\|^2 + \|\mathbf{w}_2\|^2 \leq (\|\mathbf{v}_{14}\|^2 + \|\mathbf{w}_2\|^2)^{1/2}|\mathbf{u}'\boldsymbol{\varepsilon}| + (\|\boldsymbol{\beta}_5\|_1 + \eta_1)\lambda$$
$$(5.8) \qquad\qquad + (\|\mathbf{v}_{14}\| + \|\mathbf{P}_1\mathbf{X}_2\boldsymbol{\beta}_2\|)\left(\frac{\lambda^2|A_3|}{nc_*(q_1)}\right)^{1/2},$$

where $\mathbf{u}$ is a (random) unit vector in $\mathbb{R}^n$ defined as

$$\mathbf{u} \equiv \frac{\mathbf{X}_1\boldsymbol{\Sigma}_{11}^{-1}\mathbf{Q}_{41}'\mathbf{s}_4\lambda/n - \mathbf{w}_2}{\|\mathbf{X}_1\boldsymbol{\Sigma}_{11}^{-1}\mathbf{Q}_{41}'\mathbf{s}_4\lambda/n - \mathbf{w}_2\|}. \quad (5.9)$$

Since the eigenvalues of $\boldsymbol{\Sigma}_{11}$ are no smaller than $c_*(q_1)$, we assume, without loss of generality, that $\boldsymbol{\Sigma}_{11}$ is of full rank. Since $\mathbf{X}\widehat{\boldsymbol{\beta}} = \mathbf{X}_1\widehat{\boldsymbol{\beta}}_1$ by (5.2), (5.3) gives $\mathbf{X}_1'(\mathbf{y} - \mathbf{X}_1\widehat{\boldsymbol{\beta}}_1) = \mathbf{s}_1\lambda$ so that

$$\mathbf{X}_1'\mathbf{X}_1\widehat{\boldsymbol{\beta}}_1 = \mathbf{X}_1'\mathbf{y} - \mathbf{s}_1\lambda = \mathbf{X}_1'\mathbf{X}_1\boldsymbol{\beta}_1 + \mathbf{X}_1'\mathbf{X}_2\boldsymbol{\beta}_2 + \mathbf{X}_1'\boldsymbol{\varepsilon} - \mathbf{s}_1\lambda.$$

This and the definition $\boldsymbol{\Sigma}_{jk} \equiv \mathbf{X}_j'\mathbf{X}_k/n$ yield

$$\widehat{\boldsymbol{\beta}}_1 - \boldsymbol{\beta}_1 = \boldsymbol{\Sigma}_{11}^{-1}\boldsymbol{\Sigma}_{12}\boldsymbol{\beta}_2 + \boldsymbol{\Sigma}_{11}^{-1}\mathbf{X}_1'\boldsymbol{\varepsilon}/n - \boldsymbol{\Sigma}_{11}^{-1}\mathbf{s}_1\lambda/n. \quad (5.10)$$

Inserting (5.10) into the second part of (5.1), we find that $\lambda$ is a componentwise upper bound of the absolute value of the vector

$$\mathbf{X}_2'(\mathbf{y} - \mathbf{X}\widehat{\boldsymbol{\beta}})$$
$$= \mathbf{X}_2'(\mathbf{X}_1\boldsymbol{\beta}_1 + \mathbf{X}_2\boldsymbol{\beta}_2 + \boldsymbol{\varepsilon} - \mathbf{X}_1\widehat{\boldsymbol{\beta}}_1)$$

$$= n\boldsymbol{\Sigma}_{21}\boldsymbol{\beta}_1 + n\boldsymbol{\Sigma}_{22}\boldsymbol{\beta}_2 + \mathbf{X}_2'\boldsymbol{\varepsilon}$$
$$- n\boldsymbol{\Sigma}_{21}(\boldsymbol{\beta}_1 + \boldsymbol{\Sigma}_{11}^{-1}\boldsymbol{\Sigma}_{12}\boldsymbol{\beta}_2 + \boldsymbol{\Sigma}_{11}^{-1}\mathbf{X}_1'\boldsymbol{\varepsilon}/n - \boldsymbol{\Sigma}_{11}^{-1}\mathbf{s}_1\lambda/n)$$
$$= n(\boldsymbol{\Sigma}_{22} - \boldsymbol{\Sigma}_{21}\boldsymbol{\Sigma}_{11}^{-1}\boldsymbol{\Sigma}_{12})\boldsymbol{\beta}_2 + (\mathbf{X}_2' - \boldsymbol{\Sigma}_{21}\boldsymbol{\Sigma}_{11}^{-1}\mathbf{X}_1')\boldsymbol{\varepsilon} + \boldsymbol{\Sigma}_{21}\boldsymbol{\Sigma}_{11}^{-1}\mathbf{s}_1\lambda.$$

Since $n(\boldsymbol{\Sigma}_{22} - \boldsymbol{\Sigma}_{21}\boldsymbol{\Sigma}_{11}^{-1}\boldsymbol{\Sigma}_{12}) = \mathbf{X}_2'(\mathbf{I} - \mathbf{P}_1)\mathbf{X}_2$ and $\mathbf{X}_2' - \boldsymbol{\Sigma}_{21}\boldsymbol{\Sigma}_{11}^{-1}\mathbf{X}_1' = \mathbf{X}_2'(\mathbf{I} - \mathbf{P}_1)$,

(5.11) $\quad -\lambda \le \mathbf{X}_2'(\mathbf{I} - \mathbf{P}_1)\mathbf{X}_2\boldsymbol{\beta}_2 + \mathbf{X}_2'(\mathbf{I} - \mathbf{P}_1)\boldsymbol{\varepsilon} + \boldsymbol{\Sigma}_{21}\boldsymbol{\Sigma}_{11}^{-1}\mathbf{s}_1\lambda \le \lambda.$

Taking the inner product of $\lambda\mathbf{Q}_{41}'\mathbf{s}_4$ and (5.10), we obtain, after some algebra, that, by (5.4) and Table 1,

$$\mathbf{v}_{14}'(\mathbf{v}_{13} + \mathbf{v}_{14})$$
(5.12)
$$= \mathbf{s}_4'\mathbf{Q}_{41}\boldsymbol{\Sigma}_{11}^{-1}\mathbf{s}_1\lambda^2/n$$
$$= \mathbf{s}_4'\mathbf{Q}_{41}\boldsymbol{\Sigma}_{11}^{-1}\boldsymbol{\Sigma}_{12}\boldsymbol{\beta}_2\lambda + \mathbf{s}_4'\mathbf{Q}_{41}\boldsymbol{\Sigma}_{11}^{-1}\mathbf{X}_1'\boldsymbol{\varepsilon}\lambda/n + \mathbf{s}_4'(\boldsymbol{\beta}_4 - \widehat{\boldsymbol{\beta}}_4)\lambda.$$

Similarly, the inner product of $\boldsymbol{\beta}_2$ and (5.11) yields

$$\|\mathbf{w}_2\|^2 = \boldsymbol{\beta}_2'\mathbf{X}_2'(\mathbf{I} - \mathbf{P}_1)\mathbf{X}_2\boldsymbol{\beta}_2$$
$$\le -\boldsymbol{\beta}_2'\mathbf{X}_2'(\mathbf{I} - \mathbf{P}_1)\boldsymbol{\varepsilon} - \boldsymbol{\beta}_2'\boldsymbol{\Sigma}_{21}\boldsymbol{\Sigma}_{11}^{-1}\lambda\mathbf{s}_1 + \|\boldsymbol{\beta}_2\|_1\lambda$$
$$= -\mathbf{w}_2'\boldsymbol{\varepsilon} - \mathbf{s}_1'\boldsymbol{\Sigma}_{11}^{-1}\boldsymbol{\Sigma}_{12}\boldsymbol{\beta}_2\lambda + \|\boldsymbol{\beta}_2\|_1\lambda.$$

Since $\mathbf{s}_4'\widehat{\boldsymbol{\beta}}_4 \ge 0$, by (5.1), and $\|\boldsymbol{\beta}_2\|_1 + \mathbf{s}_4'\boldsymbol{\beta}_4 \le \|\boldsymbol{\beta}_2\|_1 + \|\boldsymbol{\beta}_4\|_1 = \|\boldsymbol{\beta}_5\|_1 + \|\boldsymbol{\beta}_0\|_1 \le \|\boldsymbol{\beta}_5\|_1 + \eta_1$, by (2.4) and Table 1, the sum of (5.12) and the above inequality gives

(5.13)
$$\|\mathbf{v}_{14}\|^2 + \|\mathbf{w}_2\|^2 + \mathbf{v}_{14}'\mathbf{v}_{13}$$
$$\le (\mathbf{s}_4'\mathbf{Q}_{41}\boldsymbol{\Sigma}_{11}^{-1}\mathbf{X}_1'\lambda/n - \mathbf{w}_2')\boldsymbol{\varepsilon}$$
$$\quad - \mathbf{s}_3'\mathbf{Q}_{31}\boldsymbol{\Sigma}_{11}^{-1}\boldsymbol{\Sigma}_{12}\boldsymbol{\beta}_2\lambda + (\|\boldsymbol{\beta}_2\|_1 + \mathbf{s}_4'\boldsymbol{\beta}_4)\lambda$$
$$\le \|\mathbf{X}_1\boldsymbol{\Sigma}_{11}^{-1}\mathbf{Q}_{41}'\mathbf{s}_4\lambda/n - \mathbf{w}_2\| \cdot |\mathbf{u}'\boldsymbol{\varepsilon}|$$
$$\quad + \|\mathbf{v}_{13}\| \cdot \|\boldsymbol{\Sigma}_{11}^{-1/2}\boldsymbol{\Sigma}_{12}\boldsymbol{\beta}_2\|\sqrt{n} + (\|\boldsymbol{\beta}_5\|_1 + \eta_1)\lambda,$$

by the definition of $\mathbf{u}$ in (5.9). Since $\|\mathbf{X}_1\boldsymbol{\Sigma}_{11}^{-1/2}\mathbf{v}\|^2/n = \|\mathbf{v}\|^2$ for all $\mathbf{v} \in \mathbb{R}^{|A_1|}$ and $\mathbf{w}_2$ is orthogonal to $\mathbf{X}_1$, we find that $\|\mathbf{X}_1\boldsymbol{\Sigma}_{11}^{-1}\mathbf{Q}_{41}'\mathbf{s}_4\lambda/n - \mathbf{w}_2\| = (\|\mathbf{v}_{14}\|^2 + \|\mathbf{w}_2\|^2)^{1/2}$. Similarly, $\|\boldsymbol{\Sigma}_{11}^{-1/2}\boldsymbol{\Sigma}_{12}\boldsymbol{\beta}_2\|\sqrt{n} = \|\mathbf{P}_1\mathbf{X}_2\boldsymbol{\beta}_2\|$. Thus, by (5.13), $\|\mathbf{v}_{14}\|^2 + \|\mathbf{w}_2\|^2$ is bounded by

$$(\|\mathbf{v}_{14}\|^2 + \|\mathbf{w}_2\|^2)^{1/2}|\mathbf{u}'\boldsymbol{\varepsilon}| + (\|\boldsymbol{\beta}_5\|_1 + \eta_1)\lambda + (\|\mathbf{v}_{14}\| + \|\mathbf{P}_1\mathbf{X}_2\boldsymbol{\beta}_2\|)\|\mathbf{v}_{13}\|.$$

This implies (5.8), since, by (5.3), (5.4) and (5.6),

$$\|\mathbf{v}_{13}\|^2 = (\lambda^2/n)\mathbf{s}_3'\mathbf{Q}_{31}\boldsymbol{\Sigma}_{11}^{-1}\mathbf{Q}_{31}'\mathbf{s}_3 \le \lambda^2|A_3|/\{nc_*(q_1)\}.$$

*Step* 2. Let $B_1 \equiv (q\lambda^2/\{nc^*(q_1)\})^{1/2}$ and $B_2 \equiv (q\lambda^2/\{nc_*(q \vee q_1)\})^{1/2}$. Consider, in this step, the event

$$(5.14) \qquad |\mathbf{u}'\boldsymbol{\varepsilon}|^2 \leq \frac{\lambda^2(q_1 \vee 1)}{4nc^*(q_1)} = (q_1 \vee 1)\frac{B_1^2}{4q}.$$

We will later show that this event has high probability. We prove that, with $q_1 \equiv |A_1|$ and in the event (5.14),

$$(5.15) \qquad (q_1 - q)^+ \leq \left\{1 + 4c^*(q_1)\frac{\eta_1 n}{\lambda q} + 4\sqrt{\frac{c^*(q_1)}{c_*(q_1)}}\left(\frac{c^*(q_1)\eta_2^2 n}{\lambda^2 q}\right)^{1/2} + \frac{4c^*(q_1)}{c_*(q_1)}\right\}q,$$

provided that the $A_1$ in (5.2) contains all labels $j$ for "large" $\beta_j$,

$$(5.16) \qquad \begin{aligned}\{j: \widehat{\beta}_j(\lambda) \neq 0 \text{ or } j \notin A_0\} \\ \subseteq A_1 \subseteq \{j: |\mathbf{x}_j\{\mathbf{y} - \mathbf{X}\widehat{\boldsymbol{\beta}}(\lambda)\}| = \lambda \text{ or } j \notin A_0\}.\end{aligned}$$

Moreover, for general $A_1$ satisfying (5.2), we prove that in the event (5.14),

$$(5.17) \quad \|\mathbf{w}_2\|^2 \leq \frac{8}{3}\left(\frac{B_1^2}{4} + \eta_1\lambda + \sqrt{2}(1+\sqrt{C_5})\eta_2 B_2 + \frac{B_2^2}{2} + \frac{4}{3}C_5 B_2^2\right),$$

with $C_5 \equiv c^*(|A_5|)/c_*(|A_1 \cup A_5|)$, and, for $c_{*,5} \equiv c_*(|A_1 \cup A_5|)$,

$$(5.18) \qquad \begin{aligned} nc_{*,5}\|\boldsymbol{\beta}_5\|^2 &\leq \frac{8}{3}\left\{\frac{B_1^2}{4} + \eta_1\lambda + \eta_2\left(\frac{\lambda^2 q}{nc_*(q_1)}\right)^{1/2} + \frac{\lambda^2 q}{2nc_*(q_1)} - \frac{3\eta_2^2}{4}\right\} \\ &\quad + \left\{\frac{4}{3}\left(\frac{q\lambda^2}{nc_{*,5}}\right)^{1/2}\sqrt{1 + c^*(|A_5|)/c_*(q_1)} + 2\eta_2\right\}^2.\end{aligned}$$

By (5.14) and (5.5), we have $|\mathbf{u}'\boldsymbol{\varepsilon}|^2 \leq (\|\mathbf{v}_{14}\|^2 + B_1^2)/4$ so that

$$(\|\mathbf{v}_{14}\|^2 + \|\mathbf{w}_2\|^2)^{1/2}|\mathbf{u}'\boldsymbol{\varepsilon}| \leq \frac{1}{4}(\|\mathbf{v}_{14}\|^2 + \|\mathbf{w}_2\|^2) + |\mathbf{u}'\boldsymbol{\varepsilon}|^2$$
$$\leq \frac{1}{2}\left(\|\mathbf{v}_{14}\|^2 + \frac{\|\mathbf{w}_2\|^2 + B_1^2}{2}\right).$$

Inserting this inequality into (5.8), we find, by algebra, that

$$(5.19) \qquad \begin{aligned} \|\mathbf{v}_{14}\|^2 &+ \frac{3}{2}\|\mathbf{w}_2\|^2 \\ &\leq \frac{B_1^2}{2} + 2(\|\boldsymbol{\beta}_5\|_1 + \eta_1)\lambda + 2(\|\mathbf{v}_{14}\| + \|\mathbf{X}_2\boldsymbol{\beta}_2\|)\left(\frac{\lambda^2|A_3|}{nc_*(q_1)}\right)^{1/2}.\end{aligned}$$

We first prove (5.15) under (5.16). It follows from (5.16) and Table 1 that $A_5 = \varnothing$, so $\|\boldsymbol{\beta}_5\|_1 = 0$, $|A_3| = q \leq q_1$ and $\|\boldsymbol{\Sigma}_{11}^{-1/2}\boldsymbol{\Sigma}_{12}\boldsymbol{\beta}_2\|\sqrt{n} = \|\mathbf{P}_1\mathbf{X}_2\boldsymbol{\beta}_2\| = \|\mathbf{P}_1\mathbf{X}_6\boldsymbol{\beta}_6\| \leq \eta_2$, by (2.10). Thus, (5.19) implies

$$\|\mathbf{v}_{14}\|^2 + \frac{3}{2}\|\mathbf{w}_2\|^2 \leq \frac{B_1^2}{2} + 2\eta_1\lambda + 2(\|\mathbf{v}_{14}\| + \eta_2)B_2.$$

Since $x^2 \leq c + 2bx$ implies $x^2 \leq (b + \sqrt{b^2 + c})^2 \leq 2c + 4b^2$ for $x = \|\mathbf{v}_{14}\|$, it follows that

$$\|\mathbf{v}_{14}\|^2 \leq B_1^2 + 4\eta_1\lambda + 4\eta_2 B_2 + 4B_2^2.$$

Since $\|\mathbf{v}_{14}\|^2 \geq (q_1 - q)^+ \lambda^2/\{nc^*(q_1)\}$, by (5.5), we find, by the definition of $B_1$ and $B_2$, that

$$(q_1 - q)^+ \leq q + \frac{c^*(q_1)n}{\lambda^2}\left\{4\eta_1\lambda + 4\eta_2\left(\frac{\lambda^2 q}{c_*(q_1)n}\right)^{1/2} + \frac{4q\lambda^2}{nc_*(q_1)}\right\}.$$

This gives (5.15) by simple algebra.

For general $A_1$ satisfying (5.2), $A_5$ is no longer empty. Still, since $|A_3| + |A_5| \leq q$ by Table 1 and $\|\boldsymbol{\Sigma}_{11}^{-1/2}\boldsymbol{\Sigma}_{12}\boldsymbol{\beta}_2\|\sqrt{n} = \|\mathbf{P}_1\mathbf{X}_2\boldsymbol{\beta}_2\|$, we have, by (5.6), that

$$\left(\frac{\lambda^2|A_3|}{nc_*(q_1)}\right)^{1/2}\|\boldsymbol{\Sigma}_{11}^{-1/2}\boldsymbol{\Sigma}_{12}\boldsymbol{\beta}_2\|\sqrt{n} + \|\boldsymbol{\beta}_5\|_1\lambda$$

$$\leq \left(\frac{\lambda^2|A_3|}{nc_*(q_1)}\right)^{1/2}\|\mathbf{P}_1\mathbf{X}_2\boldsymbol{\beta}_2\| + \left(\frac{\lambda^2|A_5|}{nc_*(q)}\right)^{1/2}\|\mathbf{X}_5\boldsymbol{\beta}_5\|$$

$$\leq \left(\frac{2\lambda^2 q}{nc_*(q_1 \vee q)}\right)^{1/2}\max(\|\mathbf{P}_1\mathbf{X}_2\boldsymbol{\beta}_2\|, \|\mathbf{X}_5\boldsymbol{\beta}_5\|).$$

Moreover, it follows from Table 1, (4.1), (5.4), (5.7) and (2.10) that

$$\max(\|\mathbf{X}_2\boldsymbol{\beta}_2\|, \|\mathbf{X}_5\boldsymbol{\beta}_5\|) \leq \sqrt{nc^*(|A_5|)\|\boldsymbol{\beta}_5\|^2} + \|\mathbf{X}_6\boldsymbol{\beta}_6\|$$

$$\leq \sqrt{C_5}\|\mathbf{w}_5\| + \|\mathbf{X}_6\boldsymbol{\beta}_6\| \leq \sqrt{C_5}\|\mathbf{w}_2\| + (1 + \sqrt{C_5})\eta_2,$$

with $C_5 \equiv c^*(|A_5|)/c_*(|A_1 \cup A_5|)$. Applying these inequalities to the right-hand side of (5.19), we find that

$$\|\mathbf{v}_{14}\|^2 + \frac{3}{2}\|\mathbf{w}_2\|^2$$

$$\leq B_1^2/2 + 2\eta_1\lambda + 2\|\mathbf{v}_{14}\|\left(\frac{\lambda^2|A_3|}{nc_*(q_1)}\right)^{1/2}$$

$$+ 2(\sqrt{C_5}\|\mathbf{w}_2\| + (1 + \sqrt{C_5})\eta_2)\left(\frac{2\lambda^2 q}{nc_*(q_1 \vee q)}\right)^{1/2}$$

$$\leq B_1^2/2 + 2\eta_1\lambda + 2(1+\sqrt{C_5})\eta_2\sqrt{2}B_2$$
$$+ 2B_2(\|\mathbf{v}_{14}\| + \sqrt{2C_5}\|\mathbf{w}_2\|)$$

since $|A_3| \leq q$ and $B_2^2 \equiv \lambda^2 q/\{nc_*(q_1 \vee q)\}$. With $2\|\mathbf{v}_{14}\|B_2 \leq \|\mathbf{v}_{14}\|^2 + B_2^2$, the above inequality gives

$$\|\mathbf{w}_2\|^2 \leq (2/3)(B_1^2/2 + 2\eta_1\lambda + 2\sqrt{2}(1+\sqrt{C_5})\eta_2 B_2 + B_2^2)$$
$$+ (4/3)\sqrt{2C_5}B_2\|\mathbf{w}_2\|.$$

Since $x^2 \leq c + bx$ implies that $x^2 \leq 2c + b^2$ for $x = \|\mathbf{w}_2\|$, this gives (5.17).

The proof of (5.18) differs slightly from that of (5.17). It suffices to consider the case of $\|\boldsymbol{\beta}_5\|\sqrt{nc_{*,5}} \geq \eta_2$. By Table 1, (5.4), the definition of $\eta_2$ with (2.10) and (5.7), $\|\mathbf{w}_2\| + \eta_2 \geq \|\mathbf{w}_5\| \geq \|\boldsymbol{\beta}_5\|\sqrt{nc_{*,5}}$ with $c_{*,5} \equiv c_*(|A_1 \cup A_5|)$, so $\|\mathbf{w}_2\|^2 \geq (\|\boldsymbol{\beta}_5\|\sqrt{nc_{*,5}} - \eta_2)^2$. By (2.10) and (4.1), $\|\mathbf{X}_2\boldsymbol{\beta}_2\| \leq \eta_2 + \|\mathbf{X}_5\boldsymbol{\beta}_5\| \leq \eta_2 + \sqrt{nc^*(|A_5|)}\|\boldsymbol{\beta}_5\|$. Thus, since $2\|\mathbf{v}_{14}\|\sqrt{\lambda^2|A_3|/\{nc_*(q_1)\}} \leq \|\mathbf{v}_{14}\|^2 + \lambda^2 q/\{nc_*(q_1)\}$, (5.19) implies that

$$\frac{3}{2}(\|\boldsymbol{\beta}_5\|\sqrt{nc_{*,5}} - \eta_2)^2$$
$$\leq \frac{B_1^2}{2} + 2(\|\boldsymbol{\beta}_5\|_1 + \eta_1)\lambda + \frac{\lambda^2 q}{nc_*(q_1)}$$
$$+ 2(\eta_2 + \|\boldsymbol{\beta}_5\|\sqrt{nc^*(|A_5|)})\left(\frac{\lambda^2|A_3|}{nc_*(q_1)}\right)^{1/2}.$$

Since $\|\boldsymbol{\beta}_5\|_1^2 \leq |A_5| \cdot \|\boldsymbol{\beta}_5\|^2$ and $|A_3| + |A_5| = q$, by Cauchy–Schwarz,

$$\|\boldsymbol{\beta}_5\|_1\lambda + \|\boldsymbol{\beta}_5\|\sqrt{nc^*(|A_5|)}\left(\frac{\lambda^2|A_3|}{nc_*(q_1)}\right)^{1/2}$$
$$\leq \|\boldsymbol{\beta}_5\|\lambda(\sqrt{|A_5|} + \sqrt{c^*(|A_5|)|A_3|/c_*(q_1)})$$
$$\leq \|\boldsymbol{\beta}_5\|\lambda\sqrt{q}(1 + c^*(|A_5|)/c_*(q_1))^{1/2}.$$

It follows from the above two inequalities that

$$\|\boldsymbol{\beta}_5\|^2 nc_{*,5}$$
$$\leq \frac{2}{3}\left\{\frac{B_1^2}{2} + 2\eta_1\lambda + \frac{\lambda^2 q}{nc_*(q_1)} + 2\eta_2\left(\frac{\lambda^2 q}{nc_*(q_1)}\right)^{1/2}\right.$$
$$\left. + 2\|\boldsymbol{\beta}_5\|\lambda\sqrt{q}(1 + c^*(|A_5|)/c_*(q_1))^{1/2} + 3\eta_2\|\boldsymbol{\beta}_5\|\sqrt{nc_{*,5}}\right\} - \eta_2^2$$
$$\leq \frac{4}{3}\left(\frac{B_1^2}{4} + \eta_1\lambda + \eta_2\left(\frac{\lambda^2 q}{nc_*(q_1)}\right)^{1/2} + \frac{\lambda^2 q}{2nc_*(q_1)} - \frac{3\eta_2^2}{4}\right)$$
$$+ \|\boldsymbol{\beta}_5\|\sqrt{nc_{*,5}}\left\{\frac{4\lambda\sqrt{q}}{3\sqrt{nc_{*,5}}}(1 + c^*(|A_5|)/c_*(q_1))^{1/2} + 2\eta_2\right\}.$$

Again, since $x^2 \leq c + 2bx$ implies that $x^2 \leq 4b^2 + 2c$ for $b^2 + c \geq 0$, (5.18) follows.

*Step* 3. In this step, we find probabilistic bounds. We shall take more generous bounds $c_*(m) = c_*$ and $c^*(m) = c^*$ in (4.1) for $m \leq q^*$ with the given constants $c_*$ and $c^*$ in (2.13) and consider the event

$$(5.20) \qquad q_1 \leq |A_1 \cup A_5| \leq q^*, \qquad |\mathbf{u}'\boldsymbol{\varepsilon}|^2 \leq \frac{(q_1 \vee 1)\lambda^2}{4c^*n}.$$

In this event, we have $C_5 = C = c^*/c_*$ by (2.15) and $c_{*,5} = c_*$. Moreover, by (2.14) and the definition of $B_1$ and $B_2$ in Step 2, we have $r_1^2 = \eta_1\lambda/B_1^2$, $r_2^2 = \eta_2^2/B_1^2$ and $B_2^2 = CB_1^2$. Thus, by (2.15), (2.16) and (2.17), in the event (5.20), the assertions (5.15), (5.17) and (5.18) of Step 2 become

$$(5.21) \quad (q_1 - q)^+ + q \leq (1 + 4r_1^2 + 4\sqrt{C}r_2 + 4C)q + q = M_1^*(\lambda)q,$$

$$(5.22) \quad \|\mathbf{w}_2\|^2 \leq \frac{8}{3}\left\{\frac{1}{4} + r_1^2 + r_2\sqrt{C}(\sqrt{2} + \sqrt{2C}) + C\left(\frac{1}{2} + \frac{4}{3}C\right)\right\}B_1^2$$
$$= M_2^*(\lambda)\frac{q\lambda^2}{c^*n}$$

and

$$(5.23) \quad \begin{aligned} nc_*\|\boldsymbol{\beta}_5\|^2 &\leq \frac{8}{3}\left(\frac{1}{4} + r_1^2 + r_2\sqrt{C} + \frac{C}{2} - \frac{3r_2^2}{4}\right)B_1^2 \\ &\quad + \left(\frac{4}{3}\sqrt{C}\sqrt{1+C} + 2r_2\right)^2 B_1^2 \\ &= \frac{8}{3}\left\{\frac{1}{4} + r_1^2 + r_2\sqrt{C}(1 + 2\sqrt{1+C}) + \frac{3r_2^2}{4} + C\left(\frac{7}{6} + \frac{2}{3}C\right)\right\}B_1^2 \\ &= M_3^*(\lambda)\frac{q\lambda^2}{c^*n}. \end{aligned}$$

We note that since the constants $r_1, r_2$ and $C$ depend only on $(\lambda, q, \eta_1, \eta_2, c_*, c^*)$ and (5.16) simply requires larger $A_1$, (5.21) holds for all $A_1$ satisfying (5.2). This is not the case in Step 2 since $c_*(q_1)$ and $c^*(q_1)$ are used without (5.20). In view of (5.2), (5.4) and Table 1, (5.21), (5.22) and (5.23) match the assertions of the theorem. Thus, it remains to show that (5.20) holds for all $\lambda$ satisfying (2.19) with the probability in (2.20).

It follows from (5.9) and (5.4) that $|\mathbf{u}'\boldsymbol{\varepsilon}|$ is no greater than

$$(5.24) \quad \chi_m^* \equiv \max_{|A|=m} \max_{\mathbf{s}\in\{\pm 1\}^m} \left|\boldsymbol{\varepsilon}' \frac{\mathbf{X}_A(\mathbf{X}_A'\mathbf{X}_A)^{-1}\mathbf{s}\lambda - (\mathbf{I} - \mathbf{P}_A)\mathbf{X}\boldsymbol{\beta}}{\|\mathbf{X}_A(\mathbf{X}_A'\mathbf{X}_A)^{-1}\mathbf{s}\lambda - (\mathbf{I} - \mathbf{P}_A)\mathbf{X}\boldsymbol{\beta}\|}\right|,$$

for $q_1 = m \geq 0$. Define as Borel sets in $\mathbb{R}^{n\times(p+1)}$

$$\Omega_{m_0} \equiv \{(\mathbf{X}, \boldsymbol{\varepsilon}) : \chi_m^* \leq \sqrt{2(1+c_0)(m \vee 1)\log(p \vee a_n)} \ \forall m \geq m_0\}.$$

Since $2(1+c_0)(m \vee 1)\log(p \vee a_n) \leq (m \vee 1)\lambda^2/(4c^*n)$ by (2.19),

$$(5.25) \quad (\mathbf{X}, \boldsymbol{\varepsilon}) \in \Omega_{m_0} \quad \Rightarrow \quad |\mathbf{u}'\boldsymbol{\varepsilon}|^2 \leq \frac{(q_1 \vee 1)\lambda^2}{4c^*n} \quad \text{for } q_1 \geq m_0 \geq 0.$$

By (5.1), (5.16) and the continuity of $\widehat{\boldsymbol{\beta}}(\lambda)$ in $\lambda$, we are able to choose $A_1$ so that it changes one-at-a-time, beginning from the initial $\lambda = \infty$ with $\widehat{\boldsymbol{\beta}} = 0$ to the lower bound in (2.19). Thus, since $M_1^*(\lambda)q + 1 \leq q^*$, by (2.19) and (2.18) for such $\lambda$, and since the path of $q_1$ cannot cross the gap between $M_1^*(\lambda)q$, and $M_1^*(\lambda)q + 1$ due to the continuity of $M_1^*(\lambda)$ in $\lambda$, (5.21) and (5.25) imply that for all $\lambda$ satisfying (2.19),

$$(5.26) \quad \begin{aligned}(\mathbf{X}, \boldsymbol{\varepsilon}) &\in \Omega_q \\ &\Rightarrow \quad q_1 \equiv \#\{j : |\mathbf{x}_j(\mathbf{y} - \mathbf{X}\widehat{\boldsymbol{\beta}})| = \lambda \text{ or } j \notin A_0\} \leq M_1^*(\lambda)q.\end{aligned}$$

By (5.24), $\chi_m^*$ is the maximum of $\binom{p}{m}2^{m\vee 1}$ standard normal variables, so

$$(5.27) \quad \begin{aligned}&1 - P\{(\mathbf{X}, \boldsymbol{\varepsilon}) \in \Omega_0\} \\ &\leq \sum_{m=0}^{\infty} 2^{m\vee 1} \binom{p}{m} \exp(-(m \vee 1)(1+c_0)\log(p \vee a_n)) \\ &\leq \frac{2}{(p \vee a_n)^{1+c_0}} + \exp\left(\frac{2p}{(p \vee a_n)^{1+c_0}}\right) - 1.\end{aligned}$$

The proof is complete, since (5.20) follows from (5.25), (5.26) and (5.27). □

**6. Related results and final remarks.** In this section, we discuss some related results and make a few final remarks.

Meinshausen and Buhlmann (2006) and Zhao and Yu (2006) proved the sign-consistency $P\{\text{sgn}(\widehat{\beta}_j) = \text{sgn}(\beta_j) \ \forall j\} \to 1$, with the convention $\text{sgn}(0) \equiv 0$, for the LASSO under (2.5) and the strong irrepresentable condition

$$(6.1) \quad \|\boldsymbol{\Sigma}_{21}\boldsymbol{\Sigma}_{11}^{-1}\mathbf{s}_1\|_\infty < 1 - \kappa, \quad \text{for some } \kappa > 0,$$

where $\boldsymbol{\Sigma}_{jk} \equiv \mathbf{X}_{A_j}\mathbf{X}_{A_k}/n$ and $\mathbf{s}_1 \equiv \text{sgn}(\boldsymbol{\beta}_1)$, with $\boldsymbol{\beta}_1 \equiv (\beta_j, j \in A_1)$, $A_1 \equiv \{j : \beta_j \neq 0\}$ and $A_2 \equiv A_1^c$. We note that the definition of $A_1$ here is different from (5.2) or (5.16). Between the two papers, Zhao and Yu (2006) imposed weaker conditions on $\{n, p, q, \boldsymbol{\beta}, \lambda\}$ as

$$(6.2) \quad \lambda \geq n^{\kappa_1}\sqrt{n\log p}, \quad \min_{\beta_j \neq 0} \beta_j^2 \geq n^{\kappa_2}\frac{q\lambda^2}{n^2}, \quad n \geq n^{\kappa_3}q\log p,$$

for large $n$ and some constants $\kappa_j > 0$, where $q \equiv \#\{j : \beta_j \neq 0\}$.

Although (6.2) is not sharp, a careful study of the arguments in these two papers reveals that under (6.1), condition (6.2) can be weakened to

$$(6.3) \qquad \frac{\kappa \lambda}{\sigma} \geq \sqrt{n a_{2n}}, \mathbf{s}_1\left(\boldsymbol{\beta}_1 - \frac{\lambda}{n}\boldsymbol{\Sigma}_{11}\mathbf{s}_1\right) \geq \sigma\sqrt{a_{1n}\operatorname{diag}(\boldsymbol{\Sigma}_{11}^{-1})/n}$$

(for each component), for the sign-consistency, via (5.10) and (5.11), provided $\boldsymbol{\varepsilon} \sim N(0, \sigma^2 \mathbf{I})$, $\|\mathbf{x}_j\|^2 = n \ \forall j$, $2\log(p-q) \leq a_{2n} \to \infty$ and $2\log q \leq a_{1n} \to \infty$. This approach was taken in Wainwright (2006) under a stronger version of (6.3). Furthermore, for random designs $\mathbf{X}$ with i.i.d. Gaussian rows, Wainwright (2006) proved that the empirical version of his conditions on $\mathbf{X}$ follow from a population version of them.

Compared with these results on the sign-consistency, our focus is the properties of the model $\widehat{A}$ selected by the LASSO under milder conditions. We impose the sparse Riesz condition (2.13), instead of (6.1), to prove the rate-consistency (2.11) in Theorem 1 in terms of the sparsity, bias and the norm of missing large coefficients. We replace the $n^{\kappa_j}$, $j = 1, 2, 3$, in (6.2) by specific constants in, respectively, (2.19), Theorem 2 and Proposition 2. The second and third inequalities in (6.2) are not imposed as conditions in Theorem 1. Moreover, we allow many small nonzero coefficients, as long as the sum of their absolute values is of the order $O(q\lambda/n)$. Desirable properties of the LASSO estimator follow as in Section 3 once we establish the appropriate upper bound for the dimension $|\widehat{A}|$ of the LASSO selection.

Zhao and Yu (2006) and Zou (2006) (for fixed $p$) showed that the irrepresentable condition is necessary for the zero-consistency: $\beta_j \neq 0 \Leftrightarrow \widehat{\beta}_j \neq 0$ with high probability. It follows from the Karush–Kuhn–Tucker condition (5.1) that when $\boldsymbol{\varepsilon} = 0$, the weaker version of (6.1) with $\kappa = 0$ is necessary and sufficient for (2.2) to be zero-consistent. However, the irrepresentable condition is somewhat restrictive. As mentioned in Zhao and Yu (2006), (6.1) holds for all possible signs of $\boldsymbol{\beta}$ if and only if the norm of $\boldsymbol{\Sigma}_{21}\boldsymbol{\Sigma}_{11}^{-1}$ is less than $1 - \kappa$ as a linear mapping from $(\mathbb{R}^q, \|\cdot\|_\infty)$ to $(\mathbb{R}^{p-q}, \|\cdot\|_\infty)$. Without knowing the set $A_1$ of nonzero $\beta_j$, it is not clear how to verify (6.1), other than using simple bounds on the correlation $\mathbf{x}_j'\mathbf{x}_k$ for $j \neq k$, as in Zhao and Yu (2006). Since $\|\boldsymbol{\Sigma}_{11}^{-1}\mathbf{s}_1\|^2$ is typically of the order $\|\mathbf{s}_1\|^2 = q$, (6.1) is not a consequence of the $\ell_2$-based sparse Riesz condition (2.13) in general. For certain large data sets, it is reasonable to expect large $\|\mathbf{s}_1\|^2 = q$, even under the assumption $q \ll \min(n, p)$. In this case, (6.1) is quite restrictive.

Bunea, Tsybakov and Wegkamp (2006) and van de Geer (2007) studied convergence rates of $\|\mathbf{X}\widehat{\boldsymbol{\beta}} - \mathbf{X}\boldsymbol{\beta}\|^2$ and $\|\widehat{\boldsymbol{\beta}} - \boldsymbol{\beta}\|_1$ under the sparsity condition (2.5) and for random designs of the form $x_{ij} = \psi_j(x^i)$, where $x^i$ are i.i.d. variables and $\psi_j$ are suitable basis functions, that is, with the rows of $\mathbf{X}$ being i.i.d. copies of $(\xi_1, \ldots, \xi_p)$ as in Section 4.2. Bunea, Tsybakov and Wegkamp

(2006) obtained (3.5) and (3.6) for $\alpha = 1$ under two sets of conditions. The first set includes the lower bound $\rho_* > 0$ in (4.5), uniform upper bounds for $\|\xi_j\|_\infty$ and $q \leq c_0 \rho_* \sqrt{n/\log p}$ as in Proposition 2(ii). The second set relaxes the restriction on $q$ to $q \leq c_0 \sqrt{n/\log p}$, but relies on the correlation bound $|\mathrm{corr}(\psi_j, \psi_k)| \leq 1/(45q)$ for $\beta_k \neq 0 = \beta_j$, which has the flavor of the strong irrepresentable condition (6.1). In fact, the sample version of this condition implies $|\boldsymbol{\Sigma}_{21} \boldsymbol{\Sigma}_{11}^{-1} \mathbf{s}_1| \leq 1/\{45\lambda_{\min}(\boldsymbol{\Sigma}_{11})\}$. van de Geer (2007) considered more general forms of loss function and risk bounds under $\max_{j \leq p} \|\xi_j\|_\infty \leq K_n$. An interesting aspect of her result is the use of $D(\boldsymbol{\beta}^*)$ in place of $q$ in her version of (3.5) and (3.6), where $\boldsymbol{\beta}^*$ is the solution of (2.2) at $\mathbf{y} = \mathbf{X}\boldsymbol{\beta}$ and $D(\boldsymbol{\beta})$ is an upper bound of $(\sum_{\beta_j \neq 0} |b_j|)^2/E|\sum_j b_j \psi_j|^2$. Since $D(\boldsymbol{\beta}) = \#\{j : \beta_j \neq 0\}/\rho_*$ works under the Riesz condition and van de Geer (2007) does not assume (4.5) or (6.1), her upper bounds are indeed of a more general form than (3.5) and (3.6) when the rows of $\mathbf{X}$ are i.i.d., although the relationship of her risk bounds to $\{n, p, q\}$ is not explicit. Bounds on $\|\mathbf{X}\widehat{\boldsymbol{\beta}} - \mathbf{X}\boldsymbol{\beta}\|^2$ and $\|\widehat{\boldsymbol{\beta}} - \boldsymbol{\beta}\|_1$ do not directly imply the rate-consistency (2.11), but the converse is true for the LASSO as in Theorem 3, even for all the $\|\cdot\|_\alpha$ losses with $\alpha \geq 1$. Greenshtein and Ritov (2004) proved the persistency of a LASSO-like estimator in prediction risk under a condition on the order of $\|\boldsymbol{\beta}\|_1$ as $n \to \infty$. Since a different performance measurement is concerned, their result does not require (4.5) or (6.1).

For the estimation of $\boldsymbol{\beta}$, Donoho (2006) proved the $\ell_2$-consistency of the LASSO estimator for $p \asymp n$ when $\mathbf{X}$ is a certain normalization of a random matrix with i.i.d. $N(0,1)$ entries. Candés and Tao (2007) proved that the LASSO-like Dantzig estimator $\widetilde{\boldsymbol{\beta}}$ has the oracle property

$$\|\widetilde{\boldsymbol{\beta}} - \boldsymbol{\beta}\|^2 = O_P(1) \frac{\log p}{n} \left( \sigma^2 + \sum_{j=1}^p \beta_j^2 \wedge \sigma^2 \right)$$

under the sparsity condition (2.5) and a "uniform uncertainty principle". Since (3.6) with $\alpha = 2$ is comparable to their result, we have provided an affirmative answer to the question posed in Efron, Hastie and Tibshirani (2007), page 2363. SRC (2.13) may still hold. Recent results on random matrices are used by Candés and Tao to bound $\delta(m)$. For example, they allow $q \max_{j,k} u_{jk}^2 \asymp 1/(\log p)^4$ when $\mathbf{X}/\sqrt{n}$ is a random sample of $n$ rows from a $p \times p$ orthonormal matrix $(u_{jk})$. Their results certainly have implications on the validity of (2.13) and (4.1) for random design matrices.

Meinshausen and Yu (2006) proved that under (2.6) and certain other regularity conditions,

(6.4) $\qquad \|\widehat{\boldsymbol{\beta}} - \boldsymbol{\beta}\|^2 \leq O_P\left(\frac{\log p}{n} \frac{m_\lambda}{c_*^2(m_\lambda)}\right) + O\left(\frac{q}{m_\lambda}\right) = o_P(1),$

where $m_\lambda \equiv c^*(n \wedge p) E \|\mathbf{y}\|^2 n/\lambda^2$. They also obtained a version of (6.4), with $q/m_\lambda$ replaced by $R^2/m_\lambda^{1-p/2}$, when $c_*(m_\lambda)$ is bounded away from zero and $\boldsymbol{\beta}$ belongs to a certain weak $\ell_\alpha$-ball of radius $R$ with $0 < \alpha < 1$. In spirit, our paper and theirs both study the LASSO under conditions on the sparse eigenvalues $c_*(m)$ and $c^*(m)$, instead of (6.1), and both allow $p \gg n$ and many small nonzero coefficients. While our focus is on the properties of the selected model $\widehat{A}$ in (2.3), specifically its sparsity $|\widehat{A}|$, bias (2.8) and the norm of the missing large coefficients (2.9), theirs is on the $\ell_2$-loss $\|\widehat{\boldsymbol{\beta}} - \boldsymbol{\beta}\|^2$. Inspired by their results, and as suggested by the reviewers, we added Section 3 in the revision to discuss the implications of our results on the LASSO estimation. Still, the results in the two papers are complementary to each other. While our results are based on the upper bound (2.21) for the sparsity, Meinshausen and Yu (2006) used $|\widehat{A}| \leq c^*(|\widehat{A}|)\|\mathbf{y}\|^2 n/\lambda^2$. This is a crucial technical difference between the two papers.

Our main result asserts that as far as the rate consistency (2.11) in model selection is concerned, the performance of the LASSO for correlated designs under the sparse Riesz condition is comparable to its performance in the much simpler orthonormal designs, as in Example 1. Although the LASSO selects all coefficients of order larger than $\sqrt{q}\lambda/n$, by Theorem 2, and is sign-consistent under (6.1) and (6.3), it could miss coefficients of orders between $\sqrt{q}\lambda/n$ and the threshold level $\lambda/n$. This discrepancy with a factor of $\sqrt{q}$ is due to the interference of the estimation bias of the LASSO estimator $\widehat{\boldsymbol{\beta}}(\lambda)$ with model selection and cannot be removed for large $q$. For example, the loss measured in (2.23) cannot be recovered after the LASSO selection. A possible remedy for this discrepancy is adaptive LASSO, but for $p \gg n$ the choice of the initial estimator is unclear [Zou (2006)]. Huang, Ma and Zhang (2007) proved the sign consistency of adaptive LASSO under certain partial orthogonality condition on the pairwise correlations among vectors $\{\mathbf{y}, \mathbf{x}_1, \ldots, \mathbf{x}_p\}$. Threshold and other selection methods can be used to remove small coefficients in $\widehat{A} \cap A_0$ after LASSO selection based on the selected data $(\mathbf{y}, \mathbf{X}_{\widehat{A}})$ [cf. (3.6) for $\alpha = \infty$, Meinshausen and Yu (2006) and the references therein].

**Acknowledgments.** The authors are grateful to Nicolai Meinshausen and Bin Yu for providing an advanced copy of their paper after our paper and theirs were both presented at the Oberwolfach workshop "Qualitative Assumptions and Regularization for High-Dimensional Data" in November 2006. The authors would like to thank the Associate Editor and two referees whose comments prompted us to add Section 3 and led to several other improvements in the paper.

Department of Statistics  
Hill Center  
Busch Campus  
Rutgers University  
Piscataway, New Jersey 08854  
USA  
E-mail: czhang@stat.rutgers.edu

Department of Statistics  
and Actuarial Science  
University of Iowa  
Iowa City, Iowa 52242  
USA  
E-mail: jian@stat.uiowa.edu